\documentclass[a4]{amsart}
\usepackage{amsmath}
\usepackage{amsthm}
\usepackage{amssymb}
\usepackage{hyperref}
\usepackage{enumerate}
\usepackage{tikz}
\usetikzlibrary{cd}

\usepackage[margin=3.75cm]{geometry}

\theoremstyle{plain}
\newtheorem{theorem}{Theorem}[section]
\newtheorem{proposition}[theorem]{Proposition}
\newtheorem{corollary}[theorem]{Corollary}
\newtheorem{lemma}[theorem]{Lemma}

\theoremstyle{definition}
\newtheorem{definition}[theorem]{Definition}

\newtheorem{notconv}[theorem]{Notations and Conventions}

\theoremstyle{remark}
\newtheorem{remark}[theorem]{Remark}
\newtheorem{example}[theorem]{Example}

\newcommand{\bC}{\mathbb{C}}
\newcommand{\bP}{\mathbb{P}}
\newcommand{\bZ}{\mathbb{Z}}

\newcommand{\bI}{\mathbf{I}}
\newcommand{\bV}{\mathbf{V}}

\newcommand{\fm}{\mathfrak{m}}

\newcommand{\Iea}{{I^{\mathrm{ea}}}}
\newcommand{\Ieafix}{{I^{\mathrm{ea}}_{\mathrm{fix}}}}
\newcommand{\Iearel}[1]{{I^{\mathrm{ea}}_{#1}}}

\newcommand{\Iealog}{{I^{\mathrm{ea}}_{\mathrm{log}}}}
\newcommand{\Iealogd}[1]{{I^{\mathrm{ea}}_{\log{#1}}}}
\newcommand{\Ies}{{I^{\mathrm{es}}}}
\newcommand{\Iesfix}{{I^{\mathrm{es}}_{\mathrm{fix}}}}
\newcommand{\Iesrel}[1]{{I^{\mathrm{es}}_{#1}}}
\newcommand{\Iesrelfix}[1]{{I^{\mathrm{es}}_{\mathrm{fix}, {#1}}}}
\newcommand{\Ieslog}{{I^{\mathrm{es}}_{\mathrm{log}}}}
\newcommand{\Ieslogd}[1]{{I^{\mathrm{es}}_{\log{#1}}}}
\newcommand{\Iec}{{I^{\mathrm{ec}}}}
\newcommand{\Icd}{{I^{\mathrm{cd}}}}

\newcommand{\pair}[2]{{{#1}\supset{#2}}}
\newcommand{\tri}[3]{{{#1}\supset{#2}\supset{#3}}}

\newcommand{\taulog}{{\tau_{\mathrm{log}}}}
\newcommand{\taulogd}[1]{{\tau_{\log{#1}}}}
\newcommand{\taurel}[1]{{\tau_{#1}}}

\newcommand{\tes}{{\tau^{\mathrm{es}}}}

\newcommand{\teslog}{{\tau^{\mathrm{es}}_{\mathrm{log}}}}
\newcommand{\teslogd}[1]{{\tau^{\mathrm{es}}_{\log{#1}}}}
\newcommand{\tec}{{\tau^{\mathrm{ec}}}}

\newcommand{\mult}{{\mathrm{mt}\,}}

\newcommand{\cC}{\mathcal{C}}
\newcommand{\cD}{\mathcal{D}}

\newcommand{\cJ}{\mathcal{J}}
\newcommand{\cO}{\mathcal{O}}
\newcommand{\cP}{\mathcal{P}}
\newcommand{\cPs}{\cP_{\mathrm{s}}}

\newcommand{\cS}{\mathcal{S}}
\newcommand{\cT}{\mathcal{T}}
\newcommand{\cZ}{\mathcal{Z}}

\newcommand{\Ann}{\mathrm{Ann}}

\newcommand{\Spec}{\mathrm{Spec}\,}

\newcommand{\CatArt}{{\mathcal{A}_\mathbb{C}}}
\newcommand{\CatComp}{{\hat{\mathcal{A}}_\mathbb{C}}}
\newcommand{\CatSet}{{({\mathbf{\mathrm{Set}}})}}

\newcommand{\Deffun}{{\underline{\mathcal{D}ef}}}
\newcommand{\cAut}{\mathcal{A}ut}

\title[Planar curve singularities relative to a smooth boundary]{Planar curve singularities having constant intersection multiplicity 
with a smooth boundary}

\author{Nobuyoshi Takahashi}
\address{
Department of Mathematics, 
Graduate School of Advanced Science and Engineering, 
Hiroshima University, 
1-3-1 Kagamiyama, Higashi-Hiroshima, 
739-8526 JAPAN}
\email{tkhsnbys@hiroshima-u.ac.jp}



\begin{document}

\maketitle

\begin{abstract}
We study deformations of curve singularities in a smooth surface 
having a constant local intersection multiplicity $w$ with a smooth boundary curve. 
In particular, we consider 
the codimension of the equisingular locus 
in the semiuniversal deformation space 
and prove results on the inclusion relations between ideals 
describing different kinds of deformations. 
A classification is given of singularities expected to appear 
in codimension $\leq 3$ in a general family. 
\end{abstract}

\section{Introduction}

The deformation theory 
of reduced planar curve singularities 
is a well-developed subject 
with various applications. 
For example, the versal deformation space 
and equigeneric, equiclassical and equisingular deformations 
are reasonably well understood 
(\cite{Wahl1974}, \cite{DiazHarris1988}, \cite{GLS2007}, \cite{GLS2018}). 
As an application, 
the theory gives us the expected dimensions 
of the parameter space of curves with assigned types of singularities 
in a linear system on a surface. 
In certain cases
it has been shown that they are the actual dimensions 
 (\cite{DiazHarris1988}, \cite{GLS2018}). 

As a variant of singularity theory, 
the theory of functions and varieties relative to a boundary divisor 
has also been studied. 
Arnold considered functions relative to a smooth boundary divisor in arbitrary dimension, 
described the semiuniversal deformations and 
classified simple singularities (\cite{Arnold1978}, \cite[Ch. 17]{AGV1988}). 
For curves on surfaces, 
\cite{Wall2012} studied deformations and the determinacy problem 
relative to a possibly singular boundary curve. 

In this paper, we consider deformations of curve singularities 
in a smooth surface having a constant local intersection multiplicity $w$ 
with a smooth boundary curve. 
Our motivation comes from the enumerative geometry on log varieties. 
A pair $(X, D)$ of a variety $X$ and a reduced divisor $D$ is called a log variety, 
and can be used as a model for the ``open'' variety $X\setminus D$. 
The study of log varieties also has applications to that of usual varieties 
via degenerations. 

A natural enumerative problem is to count curves 
with a fixed tuple of tangency orders with $D$. 
For example, 
let $S$ be a smooth rational surface, 
$D\subset S$ a smooth anticanonical curve, 
$\beta$ a linear equivalence class on $S$ 
and $d:=\beta.D$. 
For a partition $d_1, \dots, d_k$ of $d$, 
we would like to count curves $C\in\beta$ of genus $0$ 
through $k$ fixed points 
such that, 
if $\nu: \overline{C}\to C$ denotes the normalization, 
$\nu^*D=\sum_{i=1}^k d_iP_i$ for some $P_i\in\overline{C}$. 

One approach to the enumeration of rational curves on a surface 
is to consider the relative compactified Jacobians 
(\cite{Beauville1990}, \cite{FGS1999}). 
In our case, 
let $\cC\to\Lambda$ be the family of integral curves, 
not equal to $D$, in the class $\beta$. 
Then the total space $\cJ(\cC/\Lambda)$ of its relative compactified Jacobian 
is a Poisson variety in a natural way (\cite{Bottacin1995b}). 
We have a natural morphism $\Lambda\to \mathrm{Sym}^d(D)$ 
given by taking the intersection. 
The fiber over a point of $\mathrm{Sym}^d(D)$ 
is the moduli space of rank-$1$ torsion-free sheaves on curves with a fixed intersection with $D$, 
and it was shown to be a symplectic variety 
(\cite{BiswasGomez2020}, \cite{Takahashi2025}). 

Thus we would like to study, 
for a fixed $\mathfrak{d}$ in the class $\beta|_D$, 
the family of curves $C$ 
in the class $\beta$ having intersections $C|_D=\mathfrak{d}$. 
In particular, we ask what singularities appear 
when $S, D$ and $\mathfrak{d}$ are general. 
The corresponding local problem is 
the study of a germ of a curve $P\in C$ in a smooth surface $S$ 
with a fixed intersection multiplicity $w$ with a smooth curve $D$, 
in the sense that $C|_D=wP$. 

Although the semiuniversal deformation of a function or a curve 
relative to a boundary 
has been already studied in \cite{Arnold1978}, \cite{AGV1988} and \cite{Wall2012}, 
here we focus on curves with a constant intersection multiplicity $w$, 
and give a more detailed description of the semiuniversal deformation. 
We take a formal coordinate system $(x, y)$ on $S$ where $D$ is defined by $y=0$ 
and $P$ by $x=y=0$. 
Then $C$ is defined by a function of the form $F=yf(x, y)+x^w$, 
and our deformation can be described 
as a perturbation of $f(x, y)$. 
An important role is played by 
the ideal 
$\Iealogd{D}(C):=\langle \partial_xF, \partial_yF, wf(x, y)-x\partial_xf(x, y)\rangle$, 
which also appeared in 
\cite{Takahashi2025} in formulating a nondegeneracy condition 
of families of curves. 
The semiuniversal deformation is given as follows. 

\begin{theorem}[Corollary \ref{cor_versal_log} with $u(x)=1$]
Let $(S, D)$ be a formal smooth surface pair 
and $C\subset S$ a reduced algebroid curve not containing $D$, 
and let $w=C.D$. 
Then the functor $\Deffun_{C\to (\tri{S}{D}{P})}^w$ 
of deformations with the constant intersection multiplicity $w$ 
has a smooth formal semiuniversal deformation space 
isomorphic to the completion at $0$ of $\cO_S/\Iealogd{D}(C)$ 
regarded as an affine space. 

If $C=\bV(F)$ with $F=yf(x, y)+x^w$, $D=\bV(y)$, 
the residue classes of $f_1, \dots, f_k\in \mathbb{C}[[x, y]]$ 
form a basis of $\mathbb{C}[[x, y]]/\Iealogd{D}(C)$ 
and $t_1, \dots, t_k$ denote formal coordinate functions, 
then a formal semiuniversal deformation is given by 
\[
y\left(f+\sum_{i=1}^k t_i f_i\right)+x^w. 
\]
\end{theorem}

Next, we study equisingular deformations of $C$ in the logarithmic context. 
They are essentially the same as the equisingular deformations of $C\cup D$, 
since equisingular deformations preserve branches. 
It is described by an ideal 
$\Ieslogd{D}(C)$ 
which is closely related to the ideal $\Iesrel{D}(C)$ 
which appeared in the proof of \cite[Proposition 2.14]{GLS2007} 
and to the equisingular ideal $\Ies(C\cup D)$. 

In \cite{DiazHarris1988}, it was shown that there are inclusion relations 
between ideals corresponding to the tangent spaces to 
equianalytic, equisingular, equiclassical and equigeneric loci. 
We have analogous inclusions between our ideals. 

\begin{theorem}[Theorem \ref{thm_inclusion}]
Let $(S, D)$ be a formal smooth surface pair 
and $C\subset S$ a reduced algebroid curve 
not containing $D$. 
Then 
\[
\Iea(C)\subseteq \Iealogd{D}(C)\subseteq \Ieslogd{D}(C)
\subseteq \Iec(C). 
\] 
\end{theorem}

In describing expected singularities in a general family, 
it is usual to fix a ``topological'' or equisingular type of a singularity 
and ask whether it appears in the family, 
and in what codimension. 
For an equisingular type given by a singularity $C$, 
it is expected that 
curves with singularities of that type appear in codimension 
$\teslogd{D}(C) := \dim_\bC \cO_S/\Ieslogd{D}(C)$. 
Theorem \ref{thm_inclusion} implies that $\teslogd{D}(C)\geq \tec(C):=\dim_\bC \cO_S/\Iec(C)$. 
In the final part, we give a classification of singularities 
with $\teslogd{D}(C)\leq 3$ using this inequality.

\medbreak
This paper is organized as follows. 
In Section 2, we recall definitions and results on deformations of planar curve singularities 
and their equisingular deformations. 
In Section 3, 
deformation functors in the log setting are studied. 
We define the deformation functors of $C$ relative to $D$ 
in a few different ways, and see how they are related. 
Ideals corresponding to tangent spaces to such deformation functors are also studied. 
In Section 4, 
equisingular deformations in the log setting are discussed. 
We prove here inclusion relations between our ideals. 
In the last section, 
singularities with $\teslogd{D}(C)\leq 3$ are classified. 

\section*{Acknowledgements}
The author is grateful to an anonymous referee 
for bringing the relevant literature on relative deformation theory 
to his attention. 
This work was supported by JSPS KAKENHI Grant Number JP22K03229.

\section{Deformations of planar curves}

In this section, we collect some definitions and results on deformations of planar 
curve singularities. 

\begin{notconv}
A \emph{planar algebroid curve} 
is a space isomorphic to the one defined by a nonzero formal power series $F\in \bC[[x, y]]$. 
We will think of it as the spectrum 
$\Spec \bC[[x, y]]/\langle F \rangle$, 
but taking the formal spectrum makes no essential difference. 
Also, $\bC$ may be replaced by any algebraically closed field 
of characteristic $0$, except when the classical topology is concerned. 

We will mostly work with formal germs of curves, surfaces, etc. at a point, 
and the closed point is omitted from the notations. 
A \emph{formal smooth curve}, resp. \emph{formal smooth surface}, will mean 
a scheme isomorphic to $\Spec \bC[[x]]$, resp. $\Spec \bC[[x, y]]$, over $\bC$. 
A \emph{formal smooth surface pair} is a pair $(S, D)$ 
where $S$ is a formal smooth surface 
and $D\subset S$ is a closed subscheme which is a formal smooth curve. 

For an affine scheme $S$, we will regard its structure sheaf $\mathcal{O}_S$ 
just as a ring, 
and refer to its element as a function. 
If $f_1, \dots, f_k$ are functions on $S$, 
$\bV(f_1, \dots, f_k)$ is 
the subscheme defined by $\langle f_1, \dots, f_k\rangle$. 
If $Z\subset S$ is a closed subscheme, 
$\bI_{Z}$ is the defining ideal. 

The symbol $\epsilon$ will always mean an element of a $\bC$-algebra 
with $\epsilon\not=0$ and $\epsilon^2=0$. 
\end{notconv}

We consider deformations over local artinian algebras and their limits 
as in \cite{Schlessinger1968}. 

\begin{definition}
Let $\CatArt$ denote the category of local artinian 
$\bC$-algebras with residue field $\bC$, 
where morphisms are local $\bC$-homomorphisms. 
Let $\CatComp$ denote the category of complete local noetherian 
$\bC$-algebras with residue field $\bC$, 
with local $\bC$-homomorphisms. 

For an object $A$ of $\CatArt$ 
and a scheme $X$ over $\bC$, 
we will denote its base change to $A$ by $X_A$. 
The closed point of the spectrum of an object of $\CatArt$ or $\CatComp$ 
will be denoted by $0$. 
For a family $\cC$ over $\Spec A$, 
the closed fiber is denoted by $\cC_0$, and similarly for morphisms. 
\end{definition}

The most basic deformation functors considered here are the following. 
We will largely follow the notation of \cite{GLS2007} and \cite{GLS2018} 
for deformation functors and related ideals. 

\begin{definition}(See \cite[II, Definition 2.1]{GLS2007}, \cite[2.2.1]{GLS2018}.)
\begin{enumerate}[(1)]
\item
Let $S$ be a formal smooth surface, 
$P$ its closed point 
and $C\subset S$ be an algebroid curve. 
If $A$ is an object of $\CatArt$, 
a \emph{deformation of $C$ \underline{in $S$}} over $A$ 
is a subscheme $\cC$ of $S_A$ flat over $A$ 
with $\cC_0 = C$. 

Let $\Deffun_{C/S}: \CatArt\to\CatSet$ and $\Deffun_{C/S}^{\textrm{fix}}: \CatArt\to\CatSet$ 
be the functors 
given by 
\begin{eqnarray*}
\Deffun_{C/S}(A) & = & \{\hbox{deformation of $C$ in $S$ over $A$}\}, \\
\Deffun_{C/S}^{\textrm{fix}}(A) & = & 
\{\hbox{deformation of $C$ in $S$ over $A$ containing $P_A$}\}. 
\end{eqnarray*}
\item
Let $C$ be an algebroid curve 
with the closed point $P$. 
If $A$ is an object of $\CatArt$, 
a \emph{deformation of $C$} over $A$ 
is a pair $(\cC, i)$ of an $A$-scheme $\cC$ flat over $A$ 
and an isomorphism $i: C\to \cC_0$ over $\bC$. 
An isomorphism of such pairs $(\cC, i)$ and $(\cC', i')$ 
is an isomorphism $\cC\to \cC'$ over $A$ compatible with $i$ and $i'$. 

A \emph{deformation of $C$ \underline{with a section}} over $A$ 
is a triple $(\cC, \cP, i)$ of an $A$-scheme $\cC$ flat over $A$, 
a closed subscheme $\cP\subset\cC$ which is a section over $\Spec A$ 
and an isomorphism $i: C\to \cC_0$ over $\bC$. 
Note that $i$ automatically maps $P$ to $\cP_0$. 
An isomorphism of such triples is additionally required to be compatible 
with the sections. 

Let $\Deffun_C: \CatArt\to\CatSet$ and 
$\Deffun_C^{\textrm{sec}}: \CatArt\to\CatSet$ be the functors 
\begin{eqnarray*}
\Deffun_C(A) & = & 
\{\hbox{deformation of $C$ over $A$}\}/\cong, \\
\Deffun_C^{\textrm{sec}}(A) & = & 
\{\hbox{deformation of $C$ with a section over $A$}\}/\cong, 
\end{eqnarray*}
where $\cong$ signifies isomorphisms in each context. 
\item
Let $\cAut_S: \CatArt\to\CatSet$ and 
$\cAut_{S\supset P}: \CatArt\to\CatSet$ be the functors 
where $\cAut_S(A)$ is the set of automorphisms 
of $S_A$ over $A$ 
which are the identity on the closed fiber, 
and 
$\cAut_{S\supset P}(A)\subset\cAut_S(A)$ consists of automorphisms 
which preserve $P_A$. 
\end{enumerate}
\end{definition}

Some of these deformation functors are not (pro-)representable, 
but they have (formal) semiuniversal deformations. 

\begin{definition}
Let $F: \CatArt\to \CatSet$ be a functor. 
We denote by $\hat{F}: \CatComp\to \CatSet$ 
its extension to $\CatComp$ given by the inverse limits 
$\hat{F}(A)=\varprojlim_I F(A/I)$ for ideals $I$ of finite colengths. 
Let $R$ be an object of $\CatComp$ and 
$\hat{\xi}_R\in \hat{F}(R)$. 
\begin{enumerate}[(1)]
\item
The functor $h_R: \CatArt\to\CatSet$ is defined by 
$h_R(A)=\mathrm{Hom}_{\CatComp}(R, A)$, 
and $\hat{\xi}_R$ is said to \emph{pro-represent} $F$ 
if the natural transformation $h_R\to F$ 
given by $h_R(A)\ni\varphi\mapsto \hat{F}(\varphi)(\hat{\xi}_R)$ 
is a natural isomorphism. 
\item
We say that $\hat{\xi}_R$ is (formally) \emph{versal} 
if the induced morphism $h_R\to F$ is smooth. 
\item
We say that $\hat{\xi}_R$ is (formally) \emph{semiuniversal} 
if it is versal and $h_R(\bC[\epsilon])\to F(\bC[\epsilon])$ 
is bijective. 
In this case, $\hat{\xi}_R$ or $R$ (or $\Spec R$) is called a \emph{hull} of $F$. 
\end{enumerate}
We say that the functor has a versal deformation, etc. 

In many cases, there is a natural extension $\tilde{F}: \CatComp\to\CatSet$ of $F$ 
and an element $\tilde{\xi}_R\in\tilde{F}(R)$ which maps to $\hat{\xi}_R$. 
In this case, we say that $\tilde{\xi}_R$ is versal, etc. 
\end{definition}

\begin{definition}
Let $F: \CatArt\to \CatSet$ be a functor 
such that $F(\bC)$ consists of one element. 
\emph{Schlessinger's conditions} (H1)--(H4) for $F$ 
are given as follows. 
\begin{enumerate}[(H1)]
\item
For a small extension $A'\to A$ 
(i.\,e. a surjective homomorphism in $\CatArt$ whose kernel is $1$-dimensional over $\bC$) 
and a homomorphism $A''\to A$ in $\CatArt$, 
the natural map 
\[
\alpha: F(A'\times_A A'')\to F(A')\times_{F(A)} F(A'') 
\]
is surjective. 
\item
If $A=\bC$ and $A'=\bC[\epsilon]$, 
$\alpha$ is bijective. 
\item
$F(\bC[\epsilon])$ is a finite dimensional $\bC$-vector space 
by the addition and the scalar multiplication 
naturally defined under the assumptions (H1) and (H2). 
\item
The map $\alpha$ is bijective 
for any small surjective map $A'\to A$ and any homomorphism $A''\to A$ in $\CatArt$. 
(Then it is bijective for any $A'\to A$ and $A''\to A$ in $\CatArt$.) 
\end{enumerate}
The functor $F$ is called a \emph{good} deformation theory 
if (H1) and (H2) are satisfied, 
and \emph{very good} 
if (H4) (and hence also (H1) and (H2)) is satisfied. 
\end{definition}

\begin{theorem}[\cite{Schlessinger1968}]\label{thm_schlessinger}
If a functor $F: \CatArt\to \CatSet$ 
satisfies Schlessinger's conditions (H1), (H2) and (H3), 
then it has a semiuniversal deformation. 
If it also satisfies (H4), it is pro-representable. 
\end{theorem}

\begin{proposition}
Let $C$ be a reduced algebroid curve in a formal smooth surface $S$. 
Then there are natural isomorphisms of functors 
$\Deffun_{C/S}/\cAut_S \cong \Deffun_C$ and 
$\Deffun_{C/S}^{\mathrm{fix}}/\cAut_{S\supset P} \cong \Deffun_C^{\mathrm{sec}}$, 
and they all have semiuniversal deformations. 
\end{proposition}
\begin{proof}
We have the isomorphism $\Deffun_{C/S}/\cAut_S \cong \Deffun_C$ 
essentially because any small deformation of $C$ is planar. 
For the second isomorphism, we use the fact that 
the section can be trivialized (\cite[Proposition II.2.2]{GLS2007}). 

The fact that $\Deffun_C$ has a semiuniversal deformation is well-known, 
and the case of $\Deffun_C^{\mathrm{sec}}$ is similar. 
See the proof of Theorem \ref{thm_versal_log_general}. 
\end{proof}

A hull of $\Deffun_C$ is called the \emph{formal semiuniversal deformation space of $C$}. 
It is smooth and can be explicitly written down. 

\begin{definition}\label{def_iea}
(See \cite[p. 30, Definition 2.2.7]{GLS2018}.) 
Let $S$ be a formal smooth surface and 
$C=\bV(F)\subset S$ a reduced algebroid curve. 
Then the \emph{equianalytic ideal}, 
or the \emph{Tjurina ideal}, of $C$ or $F$ is defined to be 
\[
\Iea(C)=\Iea(F)=\langle \partial_x F, \partial_y F, F\rangle, 
\]
where $(x, y)$ is a coordinate system on $S$. 
The \emph{fixed equianalytic/Tjurina ideal} is defined as 
\[
\Ieafix(C)=\Ieafix(F)=\langle x, y\rangle\langle\partial_x F, \partial_y F\rangle+\langle F\rangle. 
\]
These are independent of the choices of the generator $F$ and the coordinates. 

The \emph{Tjurina number} is defined to be 
\[
\tau(C)=\tau(F)=\dim_\bC \cO_S/\Iea(C). 
\]
\end{definition}

The following is well-known. 
See \cite[Theorem 1.16]{GLS2007}, for example. 
\begin{theorem}
Let $C=\bV(F)\subset S=\Spec\bC[[x, y]]$ 
be a reduced algebroid curve. 
A formal semiuniversal deformation space of $C$ is given by 
the completion at $0$ of $\mathbb{C}[[x, y]]/\Iea(F)$ 
regarded as an affine space. 
If the residue classes of $f_1, \dots, f_k\in \mathbb{C}[[x, y]]$ 
form a basis of $\mathbb{C}[[x, y]]/\Iea(F)$ 
and $t_1, \dots, t_k$ denote formal coordinate functions, 
then a formal semiuniversal deformation is given by 
\[
F+\sum_{i=1}^k t_i f_i. 
\]
\end{theorem}

\subsection{Equigeneric, equiclassical and equisingular deformations}

\begin{definition}\label{def_icd_iec}
Let $S$ be a formal smooth surface,  
$C=\bV(F)\subset S$ a reduced algebroid curve 
and $\nu: \overline{C}\to C$ its normalization. 
\begin{enumerate}[(1)]
\item
The \emph{conductor ideal} of $C$, or of $F$, is defined to be 
\[
\Icd(C)=\Icd(F)
:=
\Ann_{\cO_S}(\nu_*\cO_{\overline{C}}/ \cO_C). 
\]
The \emph{$\delta$-invariant} of $C$ is defined as 
$\delta(C)=\delta(F)=\dim_\bC \cO_S/\Icd(C)$, 
which is known to be equal to $\dim_\bC \nu_*\cO_{\overline{C}}/ \cO_C$. 
\item
The \emph{equiclassical ideal} of $C$ or $F$ is 
\[
\Iec(C)=\Iec(F)
:=
(\nu^{\#})^{-1}
(\langle \partial_x F, \partial_y F\rangle\mathcal{O}_{\overline{C}}), 
\]
where $\nu^{\#}: \cO_S\to\cO_{\overline{C}}$ is the homomorphism 
induced by $\nu$. 
We define $\tec(C)=\tec(F)=\dim_\bC\cO_S/\Iec(C)$. 
\end{enumerate}
As the notations suggest, they are independent of the choice of $F$. 
\end{definition}

For a formal family of algebroid curves, 
its equigeneric locus and equiclassical locus in the base space are defined. 
The former roughly parametrizes curves with constant $\delta(C)$ 
and the latter curves with constant $\delta(C)$ and constant ``class'' 
in the context of projective duality. 
The ideals defined above are related to these loci 
in the following way, although we do not use this in what follows. 

\begin{theorem}[{\cite[Theorem 4.15, Theorem 5.5]{DiazHarris1988}}]
Let $C\subset S$ be a planar reduced algebroid curve 
and $B$ its formal semiuniversal deformation space. 
\begin{enumerate}[(1)]
\item
The tangent cone of the equigeneric locus in $B$ is smooth, 
and under the identification of 
the tangent space to $B$ with $\cO_S/\Iea(C)$, 
it corresponds to $\Icd(C)/\Iea(C)$. 
\item
The equiclassical locus in $B$ is smooth, 
and its tangent space corresponds to $\Iec(C)/\Iea(C)$. 
\end{enumerate}
\end{theorem}

Next we recall the definition of equisingular deformations. 
Intuitively, they are topologically trivial families. 

\begin{definition}[{\cite[(3.3)]{Wahl1974}, \cite[Definition 3.13 (2)]{DiazHarris1988}, 
\cite[Definition II.2.6, II.2.7]{GLS2007}}]
Let $C\subset S$ be a reduced algebroid curve in a formal smooth surface 
and let $P$ be the closed point. 
Let $\cC$ be a deformation of $C$ over $A$ 
and $\cP\subset \cC$ a section over $A$. 
One can take a deformation $\cS$ of $S$ 
and a closed immersion $\cC\to\cS$ compatible with $C\subset S$. 
The notion of \emph{equisingularity} 
is defined inductively as follows. 
\begin{itemize}
\item
If $C$ is smooth, then $\cC$ is equisingular along $\cP$. 
\item
If $C$ is nodal, then $\cC$ is equisingular along $\cP$ 
if and only if it is equimultiple along $\cP$ over $A$, 
which in this case is equivalent to $(\cC\supset\cP)\cong (C_A\supset P_A)$ over $A$. 
\item
Otherwise, 
$\cC$ is equisingular along $\cP$ 
if and only if $\cC$ is equimultiple along $\cP$ 
and, if $\cC'$ is the \emph{reduced total transform} 
of $\cC$ on the blowup $\tilde{\cS}$ of $\cS$ along $\cP$, 
then the completion of $\cC'$  at each closed point is equisingular 
for some section lying over $\cP$. 
\end{itemize}
It turns out that $\cP$ is unique unless $C$ is smooth (see below). 
We also say $\cC$ is \emph{equisingular at $P$}, 
or just equisingular, 
without referring to the section. 

If $\cS=S_A$ and $\cP=P_A$, 
we say $\cC$ is \emph{fixed equisingular at} $P$. 

Let $\Deffun_C^{\mathrm{es}}: \CatArt\to\CatSet$ 
be the functor of isomorphism classes of equisingular deformations. 
\end{definition}

\begin{proposition}[{\cite[Theorem 3.3, Proposition 3.4]{Wahl1974}, 
\cite[Proposition II.2.66]{GLS2007}}]\label{prop_equising}
Let $\cC$ be an equisingular deformation 
of a planar reduced algebroid curve $C$ over an object $A$ of $\CatArt$, 
with an equisingular section $\cP$. 
\begin{enumerate}[(1)]
\item
The equisingular section $\cP$ is unique unless $C$ is smooth. 
\item
If $C=C_1\cup \cdots \cup C_k$ is a decomposition into nonempty unions of irreducible components, 
then there is a unique decomposition $\cC=\cC_1\cup\cdots\cup \cC_k$ 
into flat deformations $\cC_i$ of $C_i$. 
Each $\cC_i$ contains $\cP$ 
and is equisingular along $\cP$. 

(The union can be interpreted in this way: 
if $\cC$ is embedded in a deformation $\cS$ of a smooth surface, 
then $\cC_i$ is defined by an equation $\tilde{F}_i$, 
and $\cC$ is defined by $\prod_{i=1}^k \tilde{F}_i$.) 
\item
If $C_i$ is smooth and $w=C_i.C_j$($j\not=i$), 
then $\cC_i\cap \cC_j=w\cP$ as a Cartier divisor on $\cC_i$. 
\end{enumerate}
\end{proposition}

\begin{remark}
In \cite[Proposition II.2.66]{GLS2007}, 
(3) is stated when $\cC_j$ comes with an equisingular deformation of parametrization, 
but it turns out that equisingularity of equations  
and equisingularity of parametrizations 
are equivalent (\cite[Theorem II.2.64]{GLS2007}). 

One can also prove this directly by considering the blow-ups in the definition of equisingularity 
and noting that $\cC_i\cap \cC_j$ can be identified 
with the restriction of the total transform of $\cC_j$ to 
the proper transform of $\cC_i$. 
\end{remark}

\begin{theorem}[{\cite[Theorem 7.4]{Wahl1974}}]
Let $C$ be a reduced planar algebroid curve. 
Then $\Deffun_C^{\mathrm{es}}$ 
has a smooth formal semiuniversal family. 
\end{theorem}

\begin{definition}[{\cite[Definition 2.15]{GLS2007}, 
cf. \cite[Proposition 6.1]{Wahl1974}}]\label{def_ies}
Let $S$ be a formal smooth surface 
and $C=\bV(F)\subset S$ a reduced algebroid curve. 
The \emph{equisingular ideal} of $C$ is defined as 
\[
\Ies(C) = \Ies(F) := 
\{g\in \cO_S\mid \hbox{$F+\epsilon g$ is equisingular over $\bC[\epsilon]$}\}, 
\]
which is an ideal of $\cO_S$. 
It can also be characterized as the vector subspace of $\cO_S$ 
such that 
the tangent space of the equisingular locus in the semiuniversal deformation space $B$ 
corresponds to $\Ies(C)/\Iea(C)$ 
under the identification of 
the tangent space to $B$ with $\cO_S/\Iea(C)$. 

We define 
$\tes(C):=\dim_\bC \cO_S/\Ies(C)$. 

Similarly, we define the \emph{fixed equisingular ideal} of $C$ as 
\[
\Iesfix(C) = \Iesfix(F) := 
\{g\in \cO_S\mid \hbox{$F+\epsilon g$ is fixed equisingular over $\bC[\epsilon]$}\}. 
\]
\end{definition}

For the ideals introduced so far, 
the following inclusion relations hold. 

\begin{theorem}[{\cite[Corollary 4.11]{DiazHarris1988}}]\label{thm_diazharris}
For a reduced algebroid curve $C\subset S$ in a formal smooth surface, we have 
\[
\Iea(C)\subseteq \Ies(C)\subseteq \Iec(C)\subseteq \Icd(C). 
\]
\end{theorem}

The invariant $\tes(C)$ can be explicitly calculated in the following way. 

\begin{theorem}[{\cite[(II.2.8.36)]{GLS2007}, \cite[Corollary 1.1.64]{GLS2018}}]\label{thm_tes}
Let $C\subset S$ be a reduced algebroid curve 
in a formal smooth surface 
and let $q_1, \dots, q_k$ be the essential infinitely near points of $S$ 
in an embedded resolution of $C$ 
and let $m_{q_i}$ be the multiplicity of the proper transform of $C$ at $q_i$. 
Then
\[
\tes(C)=\sum_{i=1}^k\frac{m_{q_i}(m_{q_i}+1)}{2} 
- \#\{i\mid \hbox{$q_i$ is a free point}\} 
- 1. 
\]
Here, an embedded resolution means 
a resolution by blow-ups of $S$ 
such that the reduced total transform of 
$C$ is a normal crossing curve. 
An infinitely near point $q$ of $S$ on $C$ is called essential 
if the reduced total transform of $C$ is neither smooth nor normal crossing at $q$. 
It is called a free point 
if it does not lie on $2$ exceptional curves over $S$. 
\end{theorem}

\section{Deformations  of planar curves in log surfaces}

Deformation theory of functions and varieties 
relative to boundary divisors 
was studied in 
e.g. \cite{Arnold1978}, \cite[Ch. 17]{AGV1988} and \cite{Wall2012}. 
Here we focus on deformations of curve singularities in a smooth surface 
with a constant intersection multiplicity to a smooth curve. 

There are several ways to formulate the deformation problem for a curve 
in a smooth surface pair. 
We will define the corresponding deformation functors 
and show that they give equivalent descriptions of the same problem.

\begin{definition}\label{def_log_functors}
Let $(S, D)$ be a formal smooth surface pair with closed point $P$ 
and $C\subset S$ a reduced algebroid curve not containing $D$. 
Let $A$ be an object of $\CatArt$. 
\begin{enumerate}[(1)]
\item
Define $\cAut_{\tri{S}{D}{P}}(A)$ to be the set of automorphisms 
of $S_A$ over $A$ 
which are the identity on the closed fiber 
and preserve $D_A$ and $P_A$. 
Later, we will consider the functor 
\[
\Deffun_{C/S}/\cAut_{\tri{S}{D}{P}}: \CatArt \to \CatSet. 
\]
\item
A \emph{deformation of $C\to (\tri{S}{D}{P})$} over $A$ 
consists of a flat deformation $\cC$ of $C$ over $A$, 
a flat deformation $\cS$ of $S$ over $A$, 
a subscheme $\cD\subset\cS$ flat over $A$, 
a section $\cP\subset\cD$ over $A$ 
and a morphism $\tilde{\iota}: \cC\to \cS$ over $A$ 
such that the restrictions $\cD_0$, $\cP_0$ and $\tilde{\iota}_0$ to the closed fibers 
coincide with 
$D$, $P$ and the inclusion $C\hookrightarrow S$ under the given identifications. 

An isomorphism of two such deformations 
$(\cC, \cS, \cD, \cP, \tilde{\iota})$ 
and $(\cC', \cS', \cD', \cP', \tilde{\iota}')$ 
consists of isomorphisms $\alpha: \cC\overset{\sim}{\to}\cC'$ and 
$\beta: \cS\overset{\sim}{\to}\cS'$ over $A$ 
deforming the identity maps, 
commuting with $\tilde{\iota}$ and $\tilde{\iota}'$, 
such that $\beta$ sends $\cD$ and $\cP$ to $\cD'$ and $\cP'$. 
Let
\[
\Deffun_{C\to(\tri{S}{D}{P})}(A)=\{\hbox{deformation of $C\to (\tri{S}{D}{P})$ over $A$}\}/\cong. 
\]
\item
A \emph{deformation of $(\tri{C\cup D}{D}{P})$ in $(\pair{S}{P})$} over $A$ 
consists of subschemes $\cD$ and $\cZ$ of $S_A$ 
which are flat deformations of $D$ and $C\cup D$ 
satisfying $P_A\subset \cD\subset \cZ$. 
Let
\[
\Deffun_{(\tri{C\cup D}{D}{P})/(\pair{S}{P})}(A)=
\{\hbox{deformation of $(\tri{C\cup D}{D}{P})$ in $(\pair{S}{P})$ over $A$}\}. 
\]
We will consider the functor 
\[
\Deffun_{(\tri{C\cup D}{D}{P})/(\pair{S}{P})}/\cAut_{\pair{S}{P}}: \CatArt \to \CatSet. 
\]
\item
A \emph{deformation of $(\tri{C\cup D}{D}{P})\to (\pair{S}{P})$} over $A$ 
consists of flat deformations 
$\cPs\subset\cD\subset \cZ$ of 
$P\subset D\subset C\cup D$ over $A$, 
a deformation $\cP\subset\cS$ of $P\subset S$ over $A$ 
and a morphism $\tilde{\kappa}: \cZ\to \cS$ over $A$ 
such that  
$\tilde{\kappa}_0$ coincides with the inclusion $C\cup D\hookrightarrow S$ under the given identifications 
and that $\tilde{\kappa}(\cP_{\mathrm{s}})=\cP$. 

With the obvious notion of isomorphism of such deformations, 
let
\[
\Deffun_{(\tri{C\cup D}{D}{P})\to(\pair{S}{P})}(A)=
\{\hbox{deformation of $(\tri{C\cup D}{D}{P})\to (\pair{S}{P})$ over $A$}\}/\cong. 
\]
\end{enumerate}
\end{definition}

\begin{remark}
\begin{enumerate}[(1)]
\item
There are also variants without $P$. 
Here we include $P$ in the data 
since it is convenient when we consider the condition $C|_D=wP$. 
\item
These functors can be regarded as deformation functors 
of diagrams (\cite[\S 10.1]{AGV1988}, \cite[II, \S 2.3]{GLS2007}). 
\end{enumerate}
\end{remark}

\begin{proposition}\label{prop_isom_logdef}
\begin{enumerate}[(1)]
\item
The functors 
\[
\begin{array}{ll}
\Deffun_{C/S}/\cAut_{\tri{S}{D}{P}}, 
& \Deffun_{C\to (\tri{S}{D}{P})},  \\
\Deffun_{(\tri{C\cup D}{D}{P})/(\pair{S}{P})}/\cAut_{\pair{S}{P}}, 
& \Deffun_{(\tri{C\cup D}{D}{P})\to (\pair{S}{P})}  \\
\end{array}
\]
are isomorphic to each other in a natural way. 
\item
There are natural transformations, injective on each object $A$ of $\CatArt$, 
from 
these functors to 
$\Deffun_{(C\cup D)/S}^{\mathrm{fix}}/\cAut_{\pair{S}{P}}
\cong\Deffun_{C\cup D}^{\mathrm{sec}}$. 

Hence they also map to 
$\Deffun_{(C\cup D)/S}/\cAut_S\cong\Deffun_{C\cup D}$. 
\end{enumerate}
\end{proposition}
\begin{proof}
(1)
Let $A$ be an object of $\CatArt$ 
and $\fm$ its maximal ideal. 
We recall that, for a nonzero power series $F\in\bC[[x, y]]$, 
a flat deformation of $\bV(F)$ in $\Spec A[[x, y]]$ 
is given as $\bV(\tilde{F})$ for an element $\tilde{F}\in A[[x, y]]$ 
with $(\tilde{F} \mod \fm) = F$. 
Note that $\tilde{F}$ is not a zero divisor. 
Also recall that the deformation of a smooth germ $S$ is trivial. 
Further, any deformation of $(\tri{S}{D}{P})$ can be trivialized. 

Let us prove 
$\Deffun_{C/S}/\cAut_{\tri{S}{D}{P}} \cong \Deffun_{C\to (\tri{S}{D}{P})}$. 
If $\cC\subset S_A$ is a deformation of $C$ over $A$, 
let $i: \cC\to S_A$ be the inclusion morphism. 
Then $\Phi_A(\cC):=(\cC, S_A, D_A, P_A, i)$ 
is a deformation of $C\to (\tri{S}{D}{P})$ over $A$. 
If $\cC\subset S_A$ and $\cC'\subset S_A$ 
are related by an automorphism $\beta\in \cAut_{\tri{S}{D}{P}}(A)$, 
then $(\beta|_{\cC}, \beta)$ gives an isomorphism of 
$\Phi_A(\cC)$ and $\Phi_A(\cC')$. 
Conversely, if $(\alpha, \beta)$ is an isomorphism 
$\Phi_A(\cC)\overset{\sim}{\to}\Phi_A(\cC')$, 
then $\beta: S_A\to S_A$ is a family of automorphisms of $(\tri{S}{D}{P})$ over $A$ 
and $\beta(\cC)=\cC'$ holds. 
Thus $\Phi_A$ is a well-defined injective map. 

If 
$(\cC, \cS, \cD, \cP, \tilde{\iota})$ 
is a representative of an element of $ \Deffun_{C\to (\tri{S}{D}{P})}(A)$, 
then by the remark at the beginning, 
$(\cS, \cD, \cP)$ may be 
identified with $(S_A, D_A, P_A)$. 
Since $\tilde{\iota}_0$ is a closed immersion, 
so is $\tilde{\iota}$ and $\cC$ can be regarded as 
a closed subscheme of $S_A$. 
Thus $\Phi_A$ is surjective. 

Since $\Phi_A$ is obviously natural in $A$, we have 
$\Deffun_{C/S}/\cAut_{\tri{S}{D}{P}} \cong \Deffun_{C\to (\tri{S}{D}{P})}$. 

The proof of 
$\Deffun_{(\tri{C\cup D}{D}{P})/(\pair{S}{P})}/\cAut_{\pair{S}{P}}\cong\Deffun_{(\tri{C\cup D}{D}{P})\to (\pair{S}{P})}$ 
is similar. 

Finally, let us prove 
$\Deffun_{C\to (\tri{S}{D}{P})} \cong \Deffun_{(\tri{C\cup D}{D}{P})\to (\pair{S}{P})}$. 
If 
$(\cC, \cS, \cD, \cP, \tilde{\iota})$ 
represents an element of $ \Deffun_{C\to (\tri{S}{D}{P})}(A)$, 
we may think of $\cC$ as a closed subscheme of $\cS$. 
Then it is defined by an element $\tilde{F}\in \cO_{\cS}$ 
which is not a zero divisor. 
Similarly, $\cD$ is defined by an element  $\tilde{G}\in \cO_{\cS}$. 
If we set $\cZ=\bV(\tilde{F}\tilde{G})\subset \cS$ with $\tilde{\kappa}: \cZ\to\cS$ 
the inclusion morphism 
and $\cPs=\cP$, 
we have an element $(\cZ, \cD, \cPs, \cS, \cP, \tilde{\kappa})$ 
of $\Deffun_{(\tri{C\cup D}{D}{P})\to (\pair{S}{P})}(A)$. 

Conversely, if 
an element of $\Deffun_{(\tri{C\cup D}{D}{P})\to (\pair{S}{P})}(A)$ 
is given by 
$\cPs\subset\cD\subset \cZ$, $\cP\subset\cS$ and 
$\tilde{\kappa}: \cZ\to \cS$, 
identifying $\cZ$ etc. with their images, 
the ideals of $\cD$ and $\cZ$ 
are each generated by one element, say $\tilde{G}$ and $\tilde{H}$. 
Since $\cD\subset \cZ$, there exists $\tilde{F}\in \cO_{\cS}$ 
such that $\tilde{H}=\tilde{F}\tilde{G}$. 
Since $\tilde{G}$ is not a zero divisor, $\tilde{F}$ is unique 
for the chosen $\tilde{G}$ and $\tilde{H}$. 
Consequently $\cC=\bV(\tilde{F})$ does not depend 
on the choice of $\tilde{G}$ and $\tilde{H}$. 
Since $\tilde{H}$ is not a zero divisor, neither is $\tilde{F}$, 
and $\cC$ is a flat deformation of $C$. 
Thus we have an element of $ \Deffun_{C\to (\tri{S}{D}{P})}(A)$. 

These correspondences are clearly inverse to each other, 
and we have the assertion. 

\smallskip
(2)
Note that 
any deformation of $C\cup D$ over $A$ can be 
embedded into $S_A$ by choosing a lift of a basis of 
$\mathfrak{m}_{C\cup D, P}/\mathfrak{m}_{C\cup D, P}^2$ 
(if $C$ is nonempty), 
and two embeddings are related by an automorphism of $S_A$ over $A$. 

Let $(\cZ, \cD, \cPs)\to (\cS, \cP)$ represent an element of 
$\Deffun_{(\tri{C\cup D}{D}{P})\to (\pair{S}{P})}(A)$. 
Then $\cZ$ with $\cPs$ determines an element of 
$\Deffun_{C\cup D}^{\mathrm{sec}}(A)$. 
Consider two elements 
$(\cZ, \cD, \cPs)\to (\cS, \cP)$ and 
$(\cZ', \cD', \cPs')\to (\cS', \cP')$ 
of $\Deffun_{(\tri{C\cup D}{D}{P})\to (\pair{S}{P})}(A)$ 
mapping to the same element. 
By our earlier remark, we may assume that 
$\cS=\cS'=\Spec A[[x, y]]$ with $\cP=\cP'$ the origin, 
$\cZ=\cZ'$ and $\cPs=\cPs'$. 
It remains to show that $\cD=\cD'$. 
But this follows from the following fact. 
\end{proof}

\begin{lemma}\label{lem_uniq_decomposition}
Let $F, G\in \bC[[x, y]]$ be nonzero relatively prime elements, 
$A$ an object of $\CatArt$, 
$\fm$ its maximal ideal, 
and 
$\tilde{F}_1, \tilde{F}_2, \tilde{G}_1, \tilde{G}_2\in A[[x, y]]$ 
be such that $\tilde{F}_1\tilde{G}_1=\tilde{F}_2\tilde{G}_2$, 
$\tilde{F}_1\equiv\tilde{F}_2\equiv F \mod \fm$ and 
$\tilde{G}_1\equiv\tilde{G}_2\equiv G \mod \fm$. 
Then there exists a unit $\tilde{u}\in 1+\fm\cdot A[[x, y]]$ 
such that $\tilde{F}_2=\tilde{u}\tilde{F}_1$ 
and $\tilde{G}_2=\tilde{u}^{-1}\tilde{G}_1$. 
\end{lemma}
\begin{proof}
Let $\fm=I_0\supset I_1\supset\dots\supset I_n=0$ 
be a composition series, 
$A_i:=A/I_i$, 
$\fm_i:=\fm/I_i$, 
$F_k^{(i)}:=(\tilde{F}_k \mod I_i)$ 
and 
$G_k^{(i)}:=(\tilde{G}_k \mod I_i)$. 
We show that, for each $i$, 
there exists a unit $\tilde{u}_i\in 1+\fm\cdot A[[x, y]]$ 
such that $F_2^{(i)}=(\tilde{u}_i\mod I_i)F_1^{(i)}$ 
and $G_2^{(i)}=(\tilde{u}_i\mod I_i)^{-1}G_1^{(i)}$. 
For $i=0$, this is clear from the assumption. 

Assume the assertion for $i$. 
Since $I_{i+1}$ is maximal in $I_i$, 
$I_i/I_{i+1}$ is generated by an element $\delta\not=0$, 
and $\fm\delta=0$. 
By replacing $\tilde{F}_2$ and $\tilde{G}_2$ 
by $\tilde{u}_i^{-1}\tilde{F}_2$ and $\tilde{u}_i\tilde{G}_2$,  
we may assume that 
$F_2^{(i)}=F_1^{(i)}$ 
and $G_2^{(i)}=G_1^{(i)}$, 
and therefore that 
$F_2^{(i+1)}=F_1^{(i+1)}+\delta \Phi$ and 
$G_2^{(i+1)}=G_1^{(i+1)}+\delta \Gamma$ 
for some $\Phi, \Gamma\in \bC[[x, y]]$. 
Then 
\[
0=F_2^{(i+1)}G_2^{(i+1)}-F_1^{(i+1)}G_1^{(i+1)}
=\delta(G_1^{(i+1)}\Phi + F_1^{(i+1)}\Gamma)
= \delta(G\Phi + F\Gamma), 
\]
so $G\Phi+F\Gamma=0$. 
Since $F$ and $G$ are relatively prime, 
we have $\Phi=PF$ and $\Gamma=-PG$ for some $P\in \bC[[x, y]]$, 
and thus $F_2^{(i+1)}=(1+\delta P)F_1^{(i+1)}$ and 
$G_2^{(i+1)}=(1+\delta P)^{-1}G_1^{(i+1)}$. 
We have the assertion by taking a lift of $1+\delta P$ in $1+\fm\cdot A[[x, y]]$. 
\end{proof}

\begin{definition}
For a positive integer $m\leq C.D$ and $\bullet$ in 
\[
\left\{
\begin{array}{ll}
\hbox{``$C/S$''}, & \hbox{``$C\to (\tri{S}{D}{P})$''}, \\
\hbox{``$(\tri{C\cup D}{D}{P})/(\pair{S}{P})$''}, &
\hbox{``$(\tri{C\cup D}{D}{P})\to (\pair{S}{P})$''}
\end{array}
\right\}, 
\]
define $\Deffun_\bullet^m(A)$ to be the subset of $\Deffun_\bullet(A)$ 
consisting of classes of 
families with $\cC|_{D_A}\supseteq mP_A$ or $\cC|_{\cD}\supseteq m\cP$ 
in the notation of Definition \ref{def_log_functors}, 
where $P_A$ or $\cP$ is regarded as a Cartier divisor on $D_A$ or $\cD$. 
(Note that $D$ is implicit in the notation when $\bullet=\hbox{``$C/S$''}$.) 
\end{definition}

\begin{definition}
Let $\cT_{\pair{S}{D}}$ and $\cT_{\tri{S}{D}{P}}$ 
be the $\cO_S$-modules of derivations on $\cO_S$ preserving $D$ and $(D, P)$, respectively. 
\end{definition}
Explicitly, if $\cO_S=\bC[[x, y]]$, $D=\bV(y)$ and $P=\bV(x, y)$, 
then 
$\cT_{\pair{S}{D}}=\cO\partial_x + \langle y\rangle\partial_y$ 
and $\cT_{\tri{S}{D}{P}}=\langle x, y\rangle\partial_x + \langle y\rangle\partial_y$, 
where $\partial_x$ and $\partial_y$ denote the partial derivation operators. 

\begin{lemma}
If $\bV(F)=\bV(F')$, 
then 
\begin{eqnarray*}
\cT_{\pair{S}{D}} F+\langle F\rangle & = & \cT_{\pair{S}{D}} F'+\langle F'\rangle \hbox{ and} \\
\cT_{\tri{S}{D}{P}}F+\langle F\rangle & = & \cT_{\tri{S}{D}{P}}F'+\langle F'\rangle. 
\end{eqnarray*}
\end{lemma}
\begin{proof}
Straightforward. 
\end{proof}

\begin{definition}\label{def_iearel}
For $C=\bV(F)$, define 
\begin{eqnarray*}
\Iearel{D}(C) =\Iearel{D}(F) & := & \cT_{\pair{S}{D}} F+\langle F\rangle, \\
\Iearel{D, P}(C) =\Iearel{D, P}(F) & := & \cT_{\tri{S}{D}{P}}F+\langle F\rangle 
\end{eqnarray*}
and $\taurel{D}(C)=\dim_\bC \cO_S/\Iearel{D}(C)$, 
$\taurel{D, P}(C)=\dim_\bC \cO_S/\Iearel{D, P}(C)$. 
\end{definition}

\begin{theorem}\label{thm_versal_log_general}
Let $(S, D)$ be a formal smooth surface pair, 
$C=\bV(F)\subset S$ a reduced algebroid curve not containing $D$ 
and $w=C.D$. 

Then for a positive integer $m\leq w$, 
the functor $\Deffun_{C\to (\tri{S}{D}{P})}^m$ has a formally smooth semiuniversal family, 
given as follows: 
let $\bI_{mP}\subset\cO_S$ be the ideal of $mP$, as a Cartier divisor on $D$, 
and take $f_1, \dots, f_k\in \bI_{mP}$ 
whose images in 
$\bI_{mP}/\Iearel{D, P}(C)$ 
form a basis. 
Then $F+\sum_{i=1}^k t_if_i$ gives 
a formal semiuniversal family for $\Deffun_{C\to (\tri{S}{D}{P})}^m$ 
over $\Spec \bC[[t_1, \dots, t_k]]$. 
\end{theorem}
Thus, if $(x, y)$ is a coordinate system on $S$ with $D=\bV(y)$, 
the base of a semiuniversal family is smooth 
with an identification (depending on $F$) of its tangent space with 
$\langle x^m, y\rangle/\langle x\partial_x F, y\partial_xF, y\partial_yF, F\rangle$. 
\begin{proof}
We may assume that $\cO_S=\bC[[x, y]]$ 
and $D=\bV(y)$. 
Since the condition on $f_1, \dots, f_k$ 
is invariant under the multiplication by a unit, 
we may use Weierstrass preparation theorem 
and assume that $F=x^w+\sum_{i=0}^{w-1}(\sum_j c_{ij} y^j)x^i$. 
By the assumption $C.D=w$, $c_{i0}=0$ for $i=0, \dots, w-1$. 

We show that $\Deffun_{C/S}^m$ and $\cAut_{\tri{S}{D}{P}}$ are 
smooth very good deformation theories. 
It then follows by formal arguments 
that $\Deffun_{C\to (\tri{S}{D}{P})}^m\cong\Deffun_{C/S}^m/\cAut_{\tri{S}{D}{P}}$ 
is smooth and good. 

For an object $A$ of $\CatArt$, with maximal ideal $\fm$, 
an element of $\cAut_{\tri{S}{D}{P}}(A)$ is given by 
$x\mapsto x+\sum a_{ij}x^iy^j$ and $y\mapsto y(1+\sum b_{ij}x^iy^j)$ 
with $a_{ij}, b_{ij}\in\fm$ and $a_{00}=0$. 
From this it is straightforward to show that 
the natural map 
$\cAut_{\tri{S}{D}{P}}(A'\times_A A'')\to 
\cAut_{\tri{S}{D}{P}}(A')\times_{\cAut_{\tri{S}{D}{P}}(A)} \cAut_{\tri{S}{D}{P}}(A'')$ 
is bijective for any $A'\to A$ and $A''\to A$ in $\CatArt$. 
Thus $\cAut_{\tri{S}{D}{P}}$ is a smooth very good deformation functor. 

To show that $\Deffun_{C/S}^m$ is a smooth very good deformation functor, 
note the following: for any object $A$ of $\CatArt$ 
with the maximal ideal $\fm$, 
a flat deformation of $C\subset S_A$ over $A$ is defined 
by a unique element of the form 
$\tilde{F}=x^w+\sum_{i=0}^{w-1}(\sum_j (c_{ij}+a_{ij}) y^j)x^i$, 
$a_{ij}\in \fm$, 
again by Weierstrass preparation theorem (with the same $w$). 
This belongs to $\Deffun_{C/S}^m(A)$ if and only if 
$a_{i0}=0$ for $i=0, \dots, m-1$. 

Now Schlessinger's condition (H4) is easy to verify: 
if 
$\tilde{F}'=x^w+\sum_{i=0}^{w-1}(\sum_j (c_{ij}+a'_{ij}) y^j)x^i$ and 
$\tilde{F}''=x^w+\sum_{i=0}^{w-1}(\sum_j (c_{ij}+a''_{ij}) y^j)x^i$ define 
deformations over $A'$ and $A''$ which coincide over $A$, 
then $(a'_{ij}, a''_{ij})$'s belong to $A'\times_A A''$ 
and define a deformation which induces $\tilde{F}'$ and $\tilde{F}''$. 

If $A'\to A$ is a surjection in $\CatArt$ 
and $\tilde{F}=x^w+\sum_{i=0}^{w-1}(\sum_j (c_{ij}+a_{ij}) y^j)x^i$ 
defines a deformation in $\Deffun_{C/S}^m(A)$, 
then lifting $a_{ij}$ to an element $a'_{ij}$ of $A'$, 
with $a'_{i0}=0$ for $i=0, \dots, m-1$, 
we see that $\Deffun_{C/S}^m$ is smooth. 

On the other hand, since $\Deffun_{C\to (\tri{S}{D}{P})}^m(\bC[\epsilon])$ 
injects to $\Deffun_{C\cup D}^{\mathrm{sec}}(\bC[\epsilon])$, 
it is finite dimensional. 
Thus $\Deffun_{C\to (\tri{S}{D}{P})}^m$ has a formal semiuniversal family 
by Theorem \ref{thm_schlessinger}. 

Let us find the tangent space 
$\Deffun_{C/S}^m(\bC[\epsilon])/\cAut_{\tri{S}{D}{P}}(\bC[\epsilon])$. 
This can be considered as 
the space of double cosets in the following way. 
The deformation of $F$, not necessarily of Weierstrass form, 
can be written as $F+\epsilon f$ with $f\in \langle y, x^m\rangle$. 
From the left, the group of units of the form $1+\epsilon \bC[[x, y]]$ acts by multiplication. 
A unit $1+\epsilon h$ modifies $f$ to $f+hF$. 
From the right, $\cAut_{\tri{S}{D}{P}}(\bC[\epsilon])$ acts. 
Its element is given as $A\partial_x+yB\partial_y$ with $A\in \langle x, y\rangle$ 
and $B\in\bC[[x, y]]$, 
and modifies $f$ to $f+A\partial_xF+yB\partial_yF$. 
Thus our tangent space is naturally isomorphic to 
\[
\langle y, x^m\rangle /
(\langle x, y\rangle \partial_xF+\langle y\rangle \partial_y F + \langle F\rangle)
= \bI_{mP}/\Iearel{D, P}(C). 
\]

Let $f_1, \dots, f_k\in \bI_{mP}$ be functions 
giving rise to a basis of $\bI_{mP}/\Iearel{D, P}(C)$. 
Then, for any local artinian $\bC$-algebra $A$ and 
a local $\bC$-homomorphism $\varphi: R=\bC[[t_1, \dots, t_k]]\to A$, 
$F+\sum_{i=1}^k \varphi(t_i)f_i$ 
defines a deformation belonging to $\Deffun_{C/S}^m(A)$, 
and this is functorial in $A$. 
Thus we have a natural transformation $h_R\to \Deffun_{C/S}^m$. 
At $A=\bC[\epsilon]$, this is bijective. 
Since $\Deffun_{C/S}^m$ is smooth, 
we see that $R$ is a hull. 
\end{proof}

The following ideal is useful in the study of $\Deffun_{C\to (S\supset D\supset P)}^w$. 

\begin{definition}\label{def_iealog}
Let $(S, D)$ be a formal smooth surface pair 
and $C\subset S$ a reduced algebroid curve 
not containing $D$. 
Define the \emph{logarithmic equianalytic ideal} of $C$ relative to $D$ to be 
\[
\Iealogd{D}(C) 
:= \Iearel{D, P}(C): \bI_D. 
\]
If $C=\bV(F)$ and $D=\bV(y)$, 
it is also written as $\Iealogd{D}(F)$, 
and is equal to $\Iearel{D, P}(F): y$. 

We define 
$\taulogd{D}(C):=\dim_\bC \cO_S/\Iealogd{D}(C)$. 
\end{definition}

Let us prove some properties of the ideals  
defined in Definitions \ref{def_iearel} and \ref{def_iealog}. 

\begin{proposition}\label{prop_ealog}
Let $(S, D)$ be a formal smooth surface pair 
and $C\subset S$ a reduced algebroid curve 
not containing $D$. 
\begin{enumerate}[(1)]
\item
$\Iearel{D}(C)=\Iea(C\cup D): \bI_D$ 
and $\Iearel{D, P}(C) = \Ieafix(C\cup D): \bI_D$. 
\item
$\Iealogd{D}(C)=\Iearel{D}(C): \bI_D = \Iea(C\cup D): \bI_D^2 = \Ieafix(C\cup D): \bI_D^2$. 
\item
If $S=\Spec \bC[[x, y]]$, 
$C=\bV(F)$ and $D=\bV(y)$ with $F=yf(x, y)+x^wu(x)$, 
\begin{eqnarray*}
\Iealogd{D}(C) & = & 
\Iea(C)+ 
\langle wf(x, y)u(x)^{-1}-x\partial_x(f(x, y)u(x)^{-1})\rangle \\
& = & 
\langle \partial_xF, \partial_yF, wf(x, y)u(x)^{-1}-x\partial_x(f(x, y)u(x)^{-1})\rangle. 
\end{eqnarray*}
In particular, if $F=yf(x, y)+x^w$, 
\begin{eqnarray*}
\Iealogd{D}(C) & = & 
\Iea(C)+ 
\langle wf(x, y)-x\partial_xf(x, y)\rangle \\
& = & 
\langle \partial_xF, \partial_yF, wf(x, y)-x\partial_xf(x, y)\rangle. 
\end{eqnarray*}
\item
If $D=\bV(y)$, 
the map 
$m_y: \cO_S\to \cO_S: g\mapsto yg$ induces an isomorphism 
\[
\cO_S/\Iealogd{D}(C)\overset{\sim}{\to} \bI_{C\cap D}/\Iearel{D, P}(C)
\]
of $\cO_S$-modules. 
\end{enumerate}
\end{proposition}
\begin{proof}
We may assume that 
$S=\Spec \bC[[x, y]]$, $C=\bV(F)$ and $D=\bV(y)$. 

\smallskip
(1)
We have 
\begin{eqnarray*}
\Iearel{D, P}(C) & = & \langle x, y\rangle\partial_x F + \langle y\rangle \partial_y F +
  \langle F\rangle, \\
\Ieafix(C\cup D) 
  & = & \langle x, y\rangle  y\partial_x F  + \langle x, y \rangle (F+y\partial_yF) + \langle yF\rangle. 
\end{eqnarray*}
It is obvious that $y(\langle x, y\rangle \partial_x F+\langle F\rangle)\subseteq \Ieafix(C\cup D)$, 
and from $y\cdot y\partial_y F=y(F+y\partial_yF)-yF$ we see that 
$\Iearel{D, P}(C) \subseteq \Ieafix(C\cup D): y$ holds. 

Conversely, let $g$ be an element of $\Ieafix(C\cup D): y$. 
Then we can write 
\[
yg=A\partial_x(yF)+B\partial_y(yF) + C\cdot yF
= y(A\partial_x F+B\partial_yF + CF)+BF 
\]
for some $A, B, C\in\mathbb{C}[[x, y]]$ with $A, B\in\langle x, y\rangle$. 
Since $F$ is not divisible by $y$, we can write $B=yB_1$ 
and 
\[
g = A\partial_x F+B_1\cdot y\partial_yF + (C+B_1)F  \in 
\Iearel{D, P}(C). 
\]
Thus 
$\Iearel{D, P}(C) \supseteq \Ieafix(C\cup D): y$ also holds. 
The proof of 
$\Iearel{D}(C) = \Iea(C\cup D): y$
is similar. 

\smallskip
(2)
We may assume that $F=yf(x, y)+x^w$. 

By definition, 
$\Iealogd{D}(C)=\Iearel{D, P}(C): y\subseteq \Iearel{D}(C): y$. 
If $g\in \Iearel{D}(C): y$, 
we have 
\begin{eqnarray*}
yg & = & A\partial_x F+B\cdot y\partial_y F + CF \\ 
 & = & A(y\partial_x f + wx^{w-1}) + By\partial_y F + C(yf+x^w) \\ 
 & = & y(A\partial_x f + B\partial_y F + Cf) + x^{w-1}(wA + xC) 
\end{eqnarray*}
for some $A, B, C\in\mathbb{C}[[x, y]]$. 
Thus we can write $wA + xC = yA_1$. 
Then $A=w^{-1}(yA_1-xC)\in\langle x, y\rangle$, 
and $g\in \Iearel{D, P}(C): y$. 
The remaining equalities follow from (1). 

\smallskip
(3)
For the second equality, we have to show that $F$ is contained in the rightmost side. 
This follows from 
\[
y u(x)\left(wf(x, y)u(x)^{-1}-x\partial_x(f(x, y)u(x)^{-1})\right) 
+ x\partial_x F 
= (w + xu(x)^{-1}u'(x))F. 
\]
The first equality can be reduced to the case $F(x, y)=yf(x, y)+x^w$. 
From 
\[
y(wf(x, y)-x\partial_xf(x, y))
= 
w\cdot F - x\partial_x F \in \Iea(C), 
\]
we obtain 
$\Iealogd{D}(C)\supseteq
\Iea(C) + \langle wf-x\partial_xf\rangle$. 

Conversely, let $g$ be an element of $\Iealogd{D}(C)$. 
Continuing the calculation in (2), 
we have 
\begin{eqnarray*}
g & = & A\partial_x f + B\partial_yF + Cf + A_1x^{w-1} \\
& = & w^{-1}(yA_1-xC)\partial_x f + B\partial_yF + Cf + A_1x^{w-1} \\
& = & w^{-1}A_1(y\partial_x f +wx^{w-1}) + B\partial_yF + 
w^{-1}C(wf - x\partial_x f ) \\
& = &  w^{-1}A_1\partial_x F + B\partial_yF + w^{-1}C(wf - x\partial_x f ), 
\end{eqnarray*}
and $g\in \Iea(C) + \langle wf-x\partial_xf\rangle$ holds. 

\smallskip
(4)
The map $m_y$ induces 
an injective homomorphism 
$\mu_y: \cO_S/\Iealogd{D}(C)\to \cO_S/\Iearel{D, P}(C)$ 
by the definition of  $\Iealogd{D}(C)$. 
It is straightforward to see that $\Iearel{D, P}(C)\subseteq \bI_{C\cap D}$, 
and since $y\in \bI_{C\cap D}=\langle y, x^w\rangle$, 
we may think of $\mu_y$ as a map to 
$\bI_{C\cap D}/\Iearel{D, P}(C)$. 
Since $x^w\equiv -yf\mod F$, it is surjective. 
\end{proof}

\begin{remark}
Thus, $\Iealogd{D}(C)$ is the ideal which appeared in \cite[Proposition 4.4]{Takahashi2025} 
to formulate a certain nondegeneracy condition of families of curve singularities 
in relation to the relative Hilbert schemes. 
\end{remark}

We have 
the following commutative diagram 
by (1), (2) of the proposition. 
\begin{center}
\begin{tikzcd}
\cO_S/\Iealogd{D}(C) \ar[r,hook, "\times y"]\ar[rd, hook, "\times y"'] 
  & \cO_S/\Iearel{D, P}(C) \ar[r, hook, "\times y"] \ar[d, two heads] 
  & \cO_S/\Ieafix(C\cup D) \ar[d, two heads] \\
& \cO_S/\Iearel{D}(C) \ar[r, hook, "\times y"'] & \cO_S/\Iea(C\cup D) \\
\end{tikzcd}
\end{center}

\begin{corollary}\label{cor_tau_diff}
Let $(S, D)$ be a formal smooth surface pair 
with closed point $P$, 
$C\subset S$ a reduced algebroid curve 
not containing $D$ 
and $w=C.D$. 
Then 
\[
\taulogd{D}(C)=\taurel{D, P}(C)-w=\taurel{D}(C)-(w-1)=
\tau(C\cup D)-(2w-1). 
\]
\end{corollary}

\begin{proof}
By Proposition \ref{prop_ealog} (4), 
we have 
\[
\taulogd{D}(C)
=\dim_\bC \cO_S/\Iealogd{D}(C)
=\dim_\bC \cO_S/\Iearel{D, P}(C) - \dim_\bC \cO_S/\bI_{C\cap D} 
= \taurel{D, P}(C)-w. 
\]
Similarly, 
since the image of 
$\times y: \cO_S/\Iealogd{D}(C) 
\hookrightarrow
\cO_S/\Iearel{D}(C)$ 
is 
$(\langle y\rangle + \Iearel{D}(C))/\Iearel{D}(C)=
\langle y, x^{w-1} \rangle/\Iearel{D}(C)$, 
we have $\taulogd{D}(C)=\taurel{D}(C)-(w-1)$. 

The cokernel of $m_y: \cO_S/\Iearel{D}(C) \to\cO_S/\Iea(C\cup D)$ 
is $\cO_S/(\langle y\rangle + \Iea(C\cup D))$, 
and 
\[
\langle y\rangle + \Iea(C\cup D)
=
\langle y, yF, y \partial_x F, F + y\partial_y F\rangle 
=
\langle y, F\rangle 
= 
\langle y, x^w \rangle. 
\]
Thus $\taurel{D}(C)=\tau(C\cup D) - w$ 
and $\taulogd{D}(C)=\tau(C\cup D) - (2w-1)$. 
\end{proof}

\begin{corollary}\label{cor_versal_log}
Let $(S, D)$ be a formal smooth surface pair 
and $C\subset S$ a reduced algebroid curve not containing $D$, 
and let $w=C.D$. 
Then the functor $\Deffun_{C\to (\tri{S}{D}{P})}^w$ 
has a smooth formal semiuniversal deformation space 
isomorphic to the completion at $0$ of $\cO_S/\Iealogd{D}(C)$ 
regarded as an affine space. 

If $C=\bV(F)$ with $F=yf(x, y)+x^wu(x)$, $D=\bV(y)$, 
the residue classes of $f_1, \dots, f_k\in \mathbb{C}[[x, y]]$ 
form a basis of $\mathbb{C}[[x, y]]/\Iealogd{D}(C)$ 
and $t_1, \dots, t_k$ denote formal coordinate functions, 
then a formal semiuniversal deformation is given by 
\[
y\left(f+\sum_{i=1}^k t_i f_i\right)+x^wu(x). 
\]
\end{corollary}

We will refer to this family as the (formal) \emph{semiuniversal log deformation of $C$}. 

\begin{proof}
In Theorem \ref{thm_versal_log_general}, 
note that $\bI_{wP}=\bI_{C\cap D}$. 
Using the isomorphism in Proposition \ref{prop_ealog} (4), 
we see that 
\[
F+\sum_{i=1}^k t_i\cdot yf_i=y\left(f+\sum_{i=1}^k t_i f_i\right)+x^wu(x) 
\]
gives a semiuniversal family for $\Deffun_{C\to (\tri{S}{D}{P})}^w$. 
\end{proof}

\section{Equisingular deformations in the log setting}

Let us define equisingular deformation functors 
and related ideals in the log setting. 

\begin{definition}
For $\bullet=\hbox{``$C/S$''}$ or ``$(\tri{C\cup D}{D}{P})/(\pair{S}{P})$'', 
let 
$\Deffun_\bullet^{\mathrm{es, fix}}$ be the subfunctor of $\Deffun_\bullet$ 
corresponding to deformations that induce 
fixed equisingular deformations of $C\cup D$. 

For $\bullet=\hbox{``$C\to (\tri{S}{D}{P})$''}$ 
or ``$(\tri{C\cup D}{D}{P})\to (\pair{S}{P})$'', 
let 
$\Deffun_\bullet^{\mathrm{es, sec}}$ be the subfunctor of $\Deffun_\bullet$ 
corresponding to deformations that induce 
equisingular deformations of $C\cup D$ 
along the given sections. 
\end{definition}

\begin{definition}\label{def_eslog}
Let $(S, D)$ be a formal smooth surface pair 
and $C=\bV(F)\subset S$ a reduced algebroid curve not containing $D$. 
\begin{enumerate}[(1)]
\item
(\cite[proof of Proposition II.2.14]{GLS2007}) 
Let 
\[
\Iesrel{D}(C)
:=\{g\in \cO_S \mid 
\hbox{$\bV(F+\epsilon g)\cup D_{\bC[\epsilon]}\subset S_{\bC[\epsilon]}$ is equisingular}\} 
= \Ies(C\cup D): \bI_D
\]
and 
\[
\Iesrel{D, P}(C)
:=\{g\in \cO_S \mid 
\hbox{$\bV(F+\epsilon g)\cup D_{\bC[\epsilon]}\subset S_{\bC[\epsilon]}$ is fixed equisingular}\} 
= \Iesfix(C\cup D): \bI_D. 
\]
In the notation of \cite{GLS2007}, 
$\Iesrel{D, P}(C) =\Iesrelfix{D}(C)$. 
Although these ideals are given 
in a little different form in \cite[proof of Proposition II.2.14]{GLS2007}, 
the definitions coincide by \cite[Proposition II.2.66]{GLS2007}. 
\item
Let the \emph{logarithmic equisingular ideal} be 
\[
\Ieslogd{D}(C) 
=\Ieslogd{D}(F) 
:=\Iesrel{D}(C): \bI_D
=\Ies(C\cup D): \bI_D^2. 
\]

\end{enumerate}
Let 
$\teslogd{D}(C):=\dim_\bC \cO_S/\Ieslogd{D}(C)$. 
\end{definition}

\begin{proposition}\label{prop_eslog}
Let $(S, D)$ be a formal smooth surface pair 
and $C\subset S$ a reduced algebroid curve 
not containing $D$. 
Then $\Ieslogd{D}(C)=\Iesrel{D, P}(C): \bI_D=\Iesfix(C\cup D): \bI_D^2$. 
\end{proposition}

\begin{proof}
We may assume that $S=\Spec \bC[[x, y]]$, 
$D=\bV(y)$ and $C=\bV(F)$. 
We may also assume that $C$ is nonempty, 
and hence that 
$F=yf(x, y)+x^w$ with $w>0$. 

If $h\in\Ieslogd{D}(C)$, 
then by definition $y^2h\in\Ies(C\cup D)$, 
i.\,e. $\bV(y(F+\epsilon yh))$ is equisingular. 
Its equisingular section $\cP$ is contained in $\bV(y)$ 
by Proposition \ref{prop_equising} (2), 
and so is defined by $\langle y, x-\epsilon a\rangle$. 
By Proposition \ref{prop_equising} (3), 
$\langle y, F+\epsilon yh\rangle = \langle y, (x-\epsilon a)^w\rangle$. 
The left hand side is 
$\langle y, x^w\rangle$, so $a=0$. 
Thus $\bV(y(F+\epsilon yh))$ is fixed equisingular, 
and $h\in \Iesrel{D, P}(C): \bI_D=\Iesfix(C\cup D): \bI_D^2$. 
\end{proof}

\begin{proposition}
Let $(S, D)$ be a formal smooth surface pair 
and $C\subset S$ a nonempty reduced algebroid curve 
not containing $D$. 
Then the functors 
\[
\begin{array}{ll}
\Deffun_{C/S}^{\mathrm{es, fix}}/\cAut_{\tri{S}{D}{P}}, 
& \Deffun_{C\to (\tri{S}{D}{P})}^{\mathrm{es, sec}},  \\
\Deffun_{(\tri{C\cup D}{D}{P})/(\pair{S}{P})}^{\mathrm{es, fix}}/\cAut_{\pair{S}{P}}, 
& \Deffun_{(\tri{C\cup D}{D}{P})\to (\pair{S}{P})}^{\mathrm{es, sec}}  \\
\end{array}
\]
are isomorphic to each other 
and to $\Deffun_{C\cup D}^{\mathrm{es}}$ 
in a natural way. 

They can be naturally regarded as subfunctors of $\Deffun_{C\to (\tri{S}{D}{P})}^w$, 
where $w=C.D$, 
and have a smooth semiuniversal deformation space 
whose tangent space can be identified with 
$\Ieslogd{D}(C)/\Iealogd{D}(C)\cong \Ies(C\cup D)/\Iea(C\cup D)$. 
\end{proposition}
\begin{proof}
It is straightforward to see that the isomorphisms of 
Proposition \ref{prop_isom_logdef} induce 
isomorphisms between these four functors. 

Let $A$ be an object of $\CatArt$. 
Let an element of $\Deffun_{(\tri{C\cup D}{D}{P})\to (\pair{S}{P})}^{\mathrm{es, sec}}(A)$ 
be given by 
$\cPs\subset\cD\subset \cZ$, $\cP\subset\cS$ and 
$\tilde{\kappa}: \cZ\to \cS$. 
By definition, 
$\cZ$ is an equisingular deformation of $C\cup D$ along $\cPs$, 
and we have a natural map 
\[
\alpha: \Deffun_{(\tri{C\cup D}{D}{P})\to (\pair{S}{P})}^{\mathrm{es, sec}}(A)
\to \Deffun_{C\cup D}^{\mathrm{es}}(A). 
\]
Let us take another deformation 
in $\Deffun_{(\tri{C\cup D}{D}{P})\to (\pair{S}{P})}^{\mathrm{es, sec}}(A)$, 
given by 
$\cPs'\subset\cD'\subset \cZ$, $\cP'\subset\cS'$ and 
$\tilde{\kappa}': \cZ\to \cS'$ with the common $\cZ$. 
By Proposition \ref{prop_equising} (1), $\cPs=\cPs'$ holds. 
We may identify $\cS$ and $\cS'$: 
if $(x, y)$ is a formal coordinate system of $\cS$ over $A$ 
with $\cP=\bV(x, y)$, 
lifts of $x|_{\cZ}$ and $y|_{\cZ}$ to $\cS'$ 
define an isomorphism $\cS\overset{\sim}\to \cS'$ lifting $id_{\cZ}$. 
Then, by the uniqueness of decomposition of $\cZ$ 
(Lemma \ref{lem_uniq_decomposition}), 
$\cD=\cD'$ holds. 
Thus $\alpha$ is injective. 

Let $\cZ$ belong to $\Deffun_{C\cup D}^{\mathrm{es}}(A)$. 
Then $\cZ$ can be embedded into $\cS\cong \Spec A[[x, y]]$ 
since $C\cup D$ is planar.
Take the equisingular section as $\cPs$ and $\cP$. 
By Proposition \ref{prop_equising} (2) 
we have a decomposition $\cZ=\cC\cup \cD$ with $\cC$ and $\cD$ deformations 
of $C$ and $D$ with $\cPs\subset\cD$. 
Thus $\alpha$ is surjective. 

By Proposition \ref{prop_equising} (3), 
$\Deffun_{C\to (\tri{S}{D}{P})}^{\mathrm{es, sec}}$ 
is a subfunctor of 
$\Deffun_{C\to (\tri{S}{D}{P})}^w$. 
By the description of the semiuniversal log deformation in Corollary \ref{cor_versal_log}, 
we can identify the map 
$\Deffun_{C\to (\tri{S}{D}{P})}^w(\bC[\epsilon])\to 
\Deffun_{C\cup D}(\bC[\epsilon])$ 
with 
$\cO_S/\Iealogd{D}(C)\to \cO_S/\Iea(C\cup D); \overline{g}\mapsto \overline{y^2g}$, 
which is injective by Proposition \ref{prop_ealog} (2). 
Thus $\Ies(C\cup D)/\Iea(C\cup D)$, 
which corresponds to $\Deffun_{C\cup D}^{\mathrm{es}}(\bC[\epsilon])$, 
is contained in the image of this map, 
and we have an isomorphism 
$\Ieslogd{D}(C)/\Iealogd{D}(C)\cong \Ies(C\cup D)/\Iea(C\cup D)$ 
by Definition \ref{def_eslog} (2). 
\end{proof}

\begin{remark}
As in the case of the usual equisingular deformations, 
the family 
$F+y\sum_{i=1}^l t_i f_i$ is \emph{not} necessarily equisingular 
for $f_1, \dots, f_l\in\Ieslogd{D}(F)$. 
\end{remark}

Combined with Corollary \ref{cor_tau_diff}, 
this proposition implies the following. 
\begin{corollary}\label{cor_tes_diff}
Let $(S, D)$ be a formal smooth surface pair, 
$C\subset S$ a reduced algebroid curve 
not containing $D$ 
and $w=C.D$. 
Then 
\[
\teslogd{D}(C)
= \tes(C\cup D)-(2w-1). 
\]
\end{corollary}
In particular, $\teslogd{D}(C)$ is an equisingular (or topological) invariant.

\subsection{Inclusion relations between ideals}

We show that the ideal $\Ieslogd{D}(C)$ 
is contained in the equiclassical ideal $\Iec(C)$. 
This will be useful in the classification of singularities 
with a given value of $\teslogd{D}(C)$.

\begin{proposition}
In the formal smooth surface $S=\Spec \bC[[x, y]]$, 
let $D=\bV(y)$ and $P=\bV(x, y)$, 
and assume that $C=\bV(F)$ is a reduced algebroid curve 
not containing $D$. 
Let $\nu: \overline{C}\to C$ be the normalization. 
For $h\in \bC[[x, y]]$, write $\bar{h}$ for $\nu^{\#}h\in\cO_{\overline{C}}$. 
\begin{enumerate}[(1)]
\item
If $\bV(y(F+\epsilon h))$ is fixed equisingular at $P$, 
then $\bar{h}\in \Iearel{D, P}(C)\cdot\cO_{\overline{C}}$. 
\item
If $\bV(y(F+\epsilon h))$ is equisingular, 
then $\bar{h}\in \Iearel{D}(C)\cdot\cO_{\overline{C}}$. 
\item
If $C$ does not contain $D':=\bV(x)$ 
and $\bV(xy(F+\epsilon h))$ is equisingular, 
then $\bV(xy(F+\epsilon h))$ is fixed equisingular at $P$, 
and $\bar{h}\in (\cT_{\pair{S}{D\cup D'}} F)\cO_{\overline{C}}$, 
where $\cT_{\pair{S}{D\cup D'}}$ is the $\cO_S$-module of 
$\bC$-derivations of $\cO_S$ preserving $D$ and $D'$. 
\end{enumerate}
\end{proposition}
\begin{proof}
Recall that $\Iearel{D}(C)=\cT_{\pair{S}{D}}F + \langle F\rangle$ 
and $\Iearel{D, P}(C)=\cT_{\tri{S}{D}{P}}F + \langle F\rangle$. 
Therefore, 
$\Iearel{D}(C)\cdot\mathcal{O}_{\overline{C}} = 
\langle \partial_x F, y\partial_y F\rangle\mathcal{O}_{\overline{C}}$, 
$\Iearel{D, P}(C)\cdot\mathcal{O}_{\overline{C}} = 
\langle x\partial_x F, y\partial_x F, y\partial_y F\rangle\mathcal{O}_{\overline{C}}$, 
and 
$\cT_{\pair{S}{D\cup D'}} F = \langle x\partial_x F, y\partial_y F\rangle$. 
 
We will prove (1)--(3) simultaneously by an induction. 
Since the assertions are clear if $C$ is empty, 
we assume that it is not. 
For $a=1, 2, 3$, 
we say that ($a$)$_n$ holds if ($a$) holds 
when an embedded resolution of $C\cup D$ or $C\cup D\cup D'$ is obtained 
by $n$ blow-ups. 

First observe that 
(1)$_n$ implies (2)$_n$: 
if $\bV(y(F+\epsilon h))$ is equisingular, 
then by Proposition \ref{prop_equising} (2), 
the equisingular section is contained in $\bV(y)$ 
and can be written as $\bV(x-a\epsilon, y)$. 
Let $\varphi$ be the automorphism of 
$\mathbb{C}[\epsilon][[x, y]]$ given by $x\mapsto x+a\epsilon, y\mapsto y$, 
then $\varphi(F+\epsilon h)=F+\epsilon h_1$ 
with $h_1=h+a\partial_x F$ 
and it is fixed equisingular at $P$. 
If an embedded resolution of $C\cup D$ is given by $n$ blow-ups and 
if (1)$_n$ is true, then 
$\bar{h}_1\in \langle x\partial_x F, y\partial_x F, y\partial_y F\rangle\mathcal{O}_{\overline{C}}$ 
holds, 
and hence $\bar{h}\in \langle \partial_x F, y\partial_y F\rangle\mathcal{O}_{\overline{C}}$. 

For (3), we note that if $\bV(xy(F+\epsilon h))$ is equisingular along a section $\cP$, 
then $\cP$ must be contained in $\bV(x)\cap \bV(y)=P_{\bC[\epsilon]}$ 
by Proposition \ref{prop_equising} (2).  
Thus $\bV(xy(F+\epsilon h))$ is fixed equisingular at $P$. 

\smallskip
We are now ready to apply induction. 

\smallskip
(1)$_0$ 
If $k=0$, then $C\cup D$ is a node, so 
$C$ is smooth, $\partial_xF\not\in\fm_P$ and $y$ is a regular parameter on $C$. 
Thus $ \langle x\partial_x F, y\partial_x F, y\partial_y F\rangle\cO_{\overline{C}}$ 
is the maximal ideal of $\cO_{\overline{C}}$. 
The fixed equisingularity at $P$ requires that $h\in\fm_P$, 
so we have the assertion. 
By the previous remark, (2)$_0$ follows. 

\smallskip
(3)$_0$ From $\bV(F)\not=\emptyset$, 
$\bV(xyF)$ can not be nodal, and this is trivially true. 

\smallskip
Now for $n>0$, we assume that (1)$_{n-1}$, (2)$_{n-1}$ and (3)$_{n-1}$ hold. 

To prove (1)$_n$, we may apply a coordinate change $x\leadsto x+ay$, 
and assume that $\bV(x)$ is not tangent to $C$. 
Then $\bV(xy(F+\epsilon h))$ is equisingular 
and an embedded resolution of $\bV(yF)$ by $n$ blow-ups 
also gives an embedded resolution of $\bV(xyF)$ (since $n>0$). 
Thus we have only to prove (3)$_n$, 
since we already know that (1)$_n$ implies (2)$_n$. 

So, assume that $\bV(xy(F+\epsilon h))$ is fixed equisingular at $P$. 
In particular it is equimultiple at $P$, 
and if $m:=\mult F$, 
we have $\mult h\geq m$. 

Let $S_1\to S$ be the blow-up at the origin, 
$C_1$ and $D_1$ the proper transforms of $C$ and $D$ 
and $E$ the exceptional curve. 
If $C_1\cap E=\{P_1, \dots, P_k\}$, 
then $\overline{C}$ decomposes into unions $\overline{C}^{(i)}$ of smooth germs, $i=1, \dots, k$, 
with normalization maps $\nu_i$ to formal neighborhoods of $P_i$. 
We have to show that 
$(\nu_i)^{\#}h\in \langle x\partial_x F, y\partial_y F\rangle \mathcal{O}_{\overline{C}^{(i)}}$ 
for each $i$. 

By interchanging $x$ and $y$ if necessary, we may assume that $P_i$ is contained in 
the affine chart given by the coordinates $x_1, y_1$ 
with $x=x_1$ and $y=x_1y_1$. 
Here $x_1$ is a formal coordinate and $y_1$ is an affine coordinate. 
We can write $F(x_1, x_1y_1)=x_1^m F_1(x_1, y_1)$ 
and $F_1$ defines $C_1$ on this affine chart. 
Differentiating both sides with respect to $x_1$ and $y_1$, 
we obtain 
\begin{eqnarray*}
\frac{\partial F}{\partial x}(x_1, x_1y_1)
+ 
y_1\frac{\partial F}{\partial y}(x_1, x_1y_1) 
& = & 
mx_1^{m-1} F_1 + x_1^m \partial_{x_1}F_1, \\
x_1\frac{\partial F}{\partial y}(x_1, x_1y_1)
& = & 
x_1^m \partial_{y_1}F_1. 
\end{eqnarray*}
From this, it follows that 
\[
x_1^m\langle x_1\partial_{x_1} F_1, y_1\partial_{y_1} F_1\rangle \mathcal{O}_{\overline{C}^{(i)}} 
\subseteq 
\langle x\partial_x F, y\partial_y F\rangle \mathcal{O}_{\overline{C}^{(i)}}. 
\]
Let us write $h(x_1, x_1y_1)=x_1^m h_1(x_1, y_1)$. 
Then $\bV(x_1y_1 (F_1+\epsilon h_1))$ is equisingular at $P_i$ 
from the inductive definition of equisingularity. 
We show that $(\nu_i)^{\#}h_1\in
\langle x_1\partial_{x_1} F_1, y_1\partial_{y_1} F_1\rangle\mathcal{O}_{\overline{C}^{(i)}}$: 
then 
$(\nu_i)^{\#}h =(\nu_i)^{\#}(x_1^mh_1)
\in
\langle x\partial_x F, y\partial_y F\rangle \mathcal{O}_{\overline{C}^{(i)}}$ 
by the above inclusion. 

Write $P_i$ as $\bV(x_1, y_1-b)$. 
If $b=0$, 
then $P_i$ is the intersection of $D_1, E$ and $C_1$, 
and the inductive assumption (3)$_{n-1}$ 
shows that 
$(\nu_i)^{\#}h_1\in 
\langle x_1\partial_{x_1} F_1, y_1\partial_{y_1} F_1\rangle\mathcal{O}_{\overline{C}^{(i)}}$. 
If $b\not=0$, 
then $P_i\in (E\cap C_1)\setminus D_1$ 
and $\bV(x_1(F_1+\epsilon h_1))$ is equisingular at $P_i$. 
By (2)$_{n-1}$, 
$(\nu_i)^{\#}h_1$ is in $\langle \partial_{y_1}F_1, x_1\partial_{x_1}F_1\rangle\cO_{\overline{C}^{(i)}}$, 
and this ideal is equal to 
$\langle x_1\partial_{x_1} F_1, y_1\partial_{y_1} F_1\rangle\cO_{\overline{C}^{(i)}}$ 
since $y_1$ is a unit at $P_i$. 
Thus (3)$_n$ holds. 
\end{proof}

\begin{proposition}
In the formal smooth surface $S=\Spec \bC[[x, y]]$, 
let $D=\bV(y)$, $D'=\bV(x)$, $P=\bV(x, y)$ 
and let $C=\bV(F)$ be a reduced algebroid curve 
with the normalization $\nu: \overline{C}\to C$. 
Then 
$(\cT_{\pair{S}{D\cup D'}} F)\mathcal{O}_{\overline{C}}=
(y\partial_y F)\mathcal{O}_{\overline{C}}$
and 
$\Iearel{D, P}(C)\cdot\mathcal{O}_{\overline{C}}
=y \langle \partial_x F, \partial_y F\rangle\mathcal{O}_{\overline{C}}
=\bI_D\cdot \Iea(C)\mathcal{O}_{\overline{C}}$. 
\end{proposition}
\begin{proof}
For each point $R$ of $\overline{C}$ over $P$, 
we check that 
$\nu^{\#}(x\partial_x F)\in (y\partial_y F)\mathcal{O}_{\overline{C}, R}$.  
Take a local parameter $t$ at $R$ 
and let $\nu$ be given by $x=x(t)$ and $y=y(t)$ at $R$. 
We have 
\[
\frac{d}{dt} F(x(t), y(t)) = 
x'(t)\partial_x F(x(t), y(t))+y'(t)\partial_y F(x(t), y(t))=0. 
\]

If $y'(t)\partial_y F(x(t), y(t))=0$, then 
either $x'(t)=0$, hence $x(t)=0$, 
or $\partial_x F(x(t), y(t))=0$. 
In any case $\nu^{\#}(x\partial_x F)=0$ at $R$ 
and the claim is obvious. 

Otherwise, 
none of 
$x'(t)$, $y'(t)$, $\partial_x F(x(t), y(t))$ and $\partial_y F(x(t), y(t))$ 
is zero, and 
\begin{eqnarray*}
\mathrm{ord}_t\ x(t)\partial_x F(x(t), y(t))
& = & \mathrm{ord}_t\  x'(t)\partial_x F(x(t), y(t)) + 1 \\
& = & \mathrm{ord}_t\  y'(t)\partial_y F(x(t), y(t)) + 1 \\
& = & \mathrm{ord}_t\  y(t)\partial_y F(x(t), y(t)), 
\end{eqnarray*}
hence $(\cT_{\pair{S}{D\cup D'}} F)\mathcal{O}_{\overline{C}, R}=
(y\partial_y F)\mathcal{O}_{\overline{C}, R}$ 
and
$(\cT_{\pair{S}{D\cup D'}} F)\mathcal{O}_{\overline{C}}=
(y\partial_y F)\mathcal{O}_{\overline{C}}$. 

Then 
$\Iearel{D, P}(C)\cdot\mathcal{O}_{\overline{C}}
=(y\partial_x F)\mathcal{O}_{\overline{C}}  + 
(\cT_{\pair{S}{D\cup D'}} F)\mathcal{O}_{\overline{C}} 
=y \langle \partial_x F, \partial_y F\rangle\mathcal{O}_{\overline{C}}$ 
holds. 
\end{proof}

\begin{theorem}\label{thm_inclusion}
Let $(S, D)$ be a formal smooth surface pair 
and $C\subset S$ a reduced algebroid curve 
not containing $D$. 
Then 
\[
\Iea(C)\subseteq \Iealogd{D}(C)\subseteq \Ieslogd{D}(C)
\subseteq \Iec(C). 
\] 
\end{theorem}
\begin{proof}
The first inclusion follows from Proposition \ref{prop_ealog} (3) 
and the second from Proposition \ref{prop_ealog} (2) and Definition \ref{def_eslog} (2). 
We have only to show that $\Ieslogd{D}(C)\subseteq \Iec(C)$. 

Let $D=\bV(y)$ and $C=\bV(F)$ 
with $F=yf(x, y)+x^w$. 
If $h\in\Ieslogd{D}(C)$, 
then 
$y^2h\in\Iesfix(C\cup D)$ 
by Proposition \ref{prop_eslog}, 
i.\,e. $\bV(y(F+\epsilon yh))$ is fixed equisingular. 

By the previous propositions, we have 
$\overline{yh}\in y\langle \partial_xF, \partial_y F\rangle\cO_{\overline{C}}$. 
Since $\bar{y}\in\cO_{\overline{C}}$ is not a zero divisor, 
we have $\bar{h}\in \langle \partial_xF, \partial_y F\rangle\cO_{\overline{C}}$, 
i.\,e., $h\in\Iec(C)$. 
\end{proof}

\section{Singularities with $\teslogd{D}(C)\leq 3$}

Let $(S, D)$ be a formal smooth surface pair with the closed point $P$. 
In a general family of algebroid curves in $S$ 
having intersection multiplicity $w$ with $D$ at $P$, 
one expects that reduced curve singularities in an equisingular class 
appear at $P$ in codimension $\teslogd{D}(C)$, 
where $C$ is a representative of the class. 
In this section, we classify singularities with $\teslogd{D}(C)\leq 3$ 
up to analytical equivalence. 

Since we have $\tec(C)\leq\teslogd{D}(C)$ 
by Theorem \ref{thm_inclusion}, 
it suffices to consider the singularities with $\tec(C)\leq 3$. 

\begin{notconv}
\begin{enumerate}[(1)]
\item
We consider algebroid curves in $S=\Spec\bC[[x, y]]$
through $P=\bV(x, y)$, relative to the divisor $D=\bV(y)$. 
We mostly drop $D$ in $\teslogd{D}(C)$, etc. 
\item
We say that $C$ and $C'$ are \emph{equivalent} 
if there exists an automorphism of $S$ preserving $D$ (and $P$) which maps $C$ to $C'$. 
Functions $F$ and $F'$ are equivalent if $\bV(F)$ and $\bV(F')$ are. 
\item
In the calculations below, 
for an ideal $I\subseteq\bC[[x, y]]$ of finite colength, 
we will talk about a Gr\"obner basis of $I$. 
This refers to the Gr\"obner basis of $I':=I\cap \bC[x, y]$. 
The ideal $I'$ is also characterized by 
$I'\cdot\bC[[x, y]]=I$ and $\bV(I')=\{(0, 0)\}$, 
and satisfies $\bC[x, y]/I'\cong \bC[[x, y]]/I$. 
The lexicographic order on $\bC[x, y]$ with $x>y$ is denoted by $>_{x>y}$, etc. 
\end{enumerate}
\end{notconv}

\begin{lemma}\label{lem_tec_atmost_three}
Let $C$ be a reduced planar curve singularity with $\tec(C)\leq 3$. 
Then it is a simple singularity of type $A$ or $D$, 
isomorphic to one listed in Table \ref{table_tec_leq_3}. 
\begin{table}[h]
\caption{Curve singularities with $\tec(C)\leq 3$}
\label{table_tec_leq_3}
\begin{tabular}{|c|c|c|l|l|} \hline
  $\tec(C)$ &  $\delta(C)$ & Type & equation & $\Iec(C)$ \\ 
\hline \hline
 $1$ & $1$ & $A_1$ & $x^2+y^2$ & $\langle x, y\rangle$ \\
 \cline{1-5}
 $2$ & $1$ & $A_2$ & $x^3+y^2$ & $\langle x^2, y\rangle$ \\
 \cline{2-5}
        & $2$ & $A_3$ & $x^4+y^2$ & $\langle x^2, y\rangle$ \\
 \cline{1-5}
 $3$ & $2$ & $A_4$ & $x^5+y^2$ & $\langle x^3, y\rangle$ \\
 \cline{2-5}
        & $3$ & $A_5$ & $x^6+y^2$ & $\langle x^3, y\rangle$ \\
 \cline{3-5}
        &        & $D_4$ & $x^2y+y^3$ & $\langle x, y\rangle^2$  \\
         \hline
\end{tabular}
\end{table}
\end{lemma}
\begin{proof}
Since $\delta(C)\leq \tec(C)$ by Theorem \ref{thm_diazharris}, 
we have 
$\delta(C)=\sum m_q(m_q-1)/2\leq 3$, 
where $m_q$ are the multiplicities of proper transforms of $C$ 
at infinitely near singular points $q$ of $C$. 
Thus the multiplicity of $C$ is at most $3$, 
and it is easy to see that $C$ is simple of type $A_n$ with $n\leq 6$, 
of type $D_n$ with $n\leq 5$ 
or of type $E_6$. 

If $C$ is of type $A_{2k}$ or $E_6$, 
its equation is $y^2+x^{2k+1}$ or $y^3+x^4$, respectively. 
For $F=y^p+x^q$ with $\gcd(p, q)=1$, 
the equiclassical ideal $\Iec(F)$ is spanned by 
$x^iy^j$ with $pi+qj\geq pq-q$ 
by \cite[Lemma 5.7]{DiazHarris1988}. 
Then it is easy to see that 
$\tec(C)=k+1$ in the $A_{2k}$ case 
and $\tec(C)=5$ in the $E_6$ case. 

If $C$ is of type $A_{2k-1}$, we write its equation as $(y-x^k)(y+x^k)$. 
The normalization $\overline{C}$ has $2$ smooth components 
with coordinates $t_1$ and $t_2$, 
and the normalization map is given by 
$(x, y)=(t_1, t_1^k)$ and $(x, y)=(t_2, -t_2^k)$. 
The Jacobian ideal $J$ is generated by $x^{2k-1}, y$ 
and $J\cdot\cO_{\overline{C}}=\langle t_1^k\rangle \times \langle t_2^k\rangle$. 
Thus $\Iec(C)=\cO_C\cap J\cdot\cO_{\overline{C}}$ 
is generated by $y, x^k$. 

If $C$ is of type $D_4$ with the equation $xy(x+y)$, 
the normalization $\overline{C}$ has $3$ components, 
with coordinates $t_1, t_2, t_3$. 
The Jacobian ideal $J$ is generated by $xy, x^2+3y^2$ 
and 
$J\cdot\cO_{\overline{C}}=
\langle t_1^2\rangle \times \langle t_2^2\rangle\times \langle t_3^2\rangle$. 
Since the branches of $C$ have distinct tangent directions, 
$\Iec(C)=\cO_C\cap J\cdot\cO_{\overline{C}}=\langle x, y\rangle^2$. 

If $C$ is of type $D_5$ with the equation $x^2y+y^4=y(x^2+y^3)$, 
then $\overline{C}$ has $2$ components 
and the normalization map is given by 
$(x, y)=(0, t_1)$ and $(x, y)=(t_2^3, -t_2^2)$. 
The Jacobian ideal is generated by $xy, x^2+4y^3$ 
and $J\cdot\cO_{\overline{C}}=\langle t_1^2\rangle \times \langle t_2^5\rangle$. 
We see that 
$\Iec(C)=\langle x^2, xy, y^3\rangle$ and $\tec(C)=4$. 
\end{proof}

For later use, we calculate the equianalytic invariants 
of certain polynomials. 
\begin{lemma}\label{lem_binom}
Let $F=y^l+x^k$. 
Then $\Iea(F)=\Iealog(F)=\langle y^{l-1}, x^{k-1}\rangle$ 
and $\tau(F)=\taulog(F)=(k-1)(l-1)$. 
\end{lemma}
\begin{proof}
Immediate from Proposition \ref{prop_ealog} (3). 
\end{proof}

\begin{lemma}\label{lem_trinom}
For $k>k'>0$ and $l\geq l'\geq l/2$ 
satisfying $kl'+k'l\not=kl$, 
and for any $a\in\bC\setminus \{0\}$, 
$y^l+ay^{l'}x^{k'}+x^k$ is equivalent to 
$F=y^l+y^{l'}x^{k'}+x^k$. 

A Gr\"obner basis of $\Iealog(F)$ for $>_{x>y}$ is given by 
$\{kx^{k-1}+k'x^{k'-1}y^{l'}, x^{k'}y^{l'-1}, y^{l-1}\}$, 
and $\taulog(F)=(k-1)(l'-1)+k'(l-l')$. 
\end{lemma}
\begin{proof}
The first statement follows from a scaling of $x, y$. 

For 
$f_1:=\partial_xF=kx^{k-1}+k'x^{k'-1}y^{l'}$, 
$f_2:=\partial_yF=l'x^{k'}y^{l'-1}+ly^{l-1}$ and 
$f_3:=k(y^{l-1}+y^{l'-1}x^{k'})-x\partial_x(y^{l-1}+y^{l'-1}x^{k'})
= (k-k')x^{k'}y^{l'-1}+ky^{l-1}$, 
we have $\langle f_2, f_3\rangle = \langle x^{k'}y^{l'-1}, y^{l-1}\rangle$ 
by the assumption $kl'+k'l\not=kl$. 
Thus 
$\Iealog(F)=\langle f_1, x^{k'}y^{l'-1}, y^{l-1}\rangle$. 
The $S$-polynomial of $f_1$ and $x^{k'}y^{l'-1}$ is, up to a scalar, 
$y^{l'-1}f_1-kx^{k-k'-1}\cdot x^{k'}y^{l'-1} = k'x^{k'-1}y^{2l'-1}$, 
which is a multiple of $y^{l-1}$ by the assumption $l'\geq l/2$. 
It is easy to see that $S(f_1, y^{l-1})\in\langle y^{l-1}\rangle$, 
so $\{kx^{k-1}+k'x^{k'-1}y^{l'}, x^{k'}y^{l'-1}, y^{l-1}\}$ 
is a Gr\"obner basis 
(though not necessarily minimal). 

The monomials $x^iy^j$ in the following range form 
a basis of $\bC[[x, y]]/\Iealog(F)$: 
$0\leq j<l'-1$ and $0\leq i<k-1$, 
or $l'-1\leq j<l-1$ and $0\leq i<k'$. 
\end{proof}

The following is useful in simplifying the equations.  
\begin{lemma}\label{lem_scaling}
For $U(x)\in\bC[[x]]^\times$, 
$d\in\bZ\setminus\{0\}$ and $e\in\bZ$, 
there exists $v(x)\in\bC[[x]]^\times$ 
such that 
$v(x)^dU(xv(x)^e)=1$. 
\end{lemma}
\begin{proof}
Let $U(x)=\sum_{i=0}^\infty a_ix^i$. 
We will inductively find $s_i\in\bC$ such that, 
if we set $v_k(x)=\sum_{i=0}^k s_ix^i$, 
$v_k(x)^dU(xv_k(x)^e)\equiv 1\mod x^{k+1}$. 

For $k=0$, this amounts to taking an $s_0$ such that $s_0^d=a_0^{-1}$. 

For $k>0$, assume that $s_0, \dots, s_{k-1}$ have been chosen so that 
$v_{k-1}(x)^dU(xv_{k-1}(x)^e)\equiv 1\mod x^k$. 
Then we have $v_{k-1}(x)^dU(xv_{k-1}(x)^e)\equiv 1+cx^k\mod x^{k+1}$ 
for some $c\in\bC$. 
We have 
\begin{eqnarray*}
v_k(x)^d U(xv_k(x)^e) & = & 
(v_{k-1}(x)+s_kx^k)^dU(x(v_{k-1}(x)+s_kx^k)^e) \\
& \equiv & 
v_{k-1}(x)^dU(xv_{k-1}(x)^e)+dv_{k-1}(x)^{d-1}s_kx^kU(xv_{k-1}(x)^e) \mod x^{k+1} \\
& \equiv & 
1+cx^k +dv_{k-1}(0)^{d-1}U(0)s_kx^k  \mod x^{k+1} \\
& = & 
1 + (c+ds_0^{-1}s_k)x^k. 
\end{eqnarray*}
Thus it suffices to take $s_k$ to be $-d^{-1}s_0c$. 
\end{proof}

The following is a typical singularity in the log setting. 
\begin{proposition}\label{prop_two_branches}

Let $C$ be an $A_{2k-1}$-singularity, isomorphic to $\bV(y^2+x^{2k})$, 
and assume that $w := C.D > 2k$. 

Then $C$ is equivalent to $\bV(F)$ with $F=(y+x^k)(y+x^{w-k})=y(y+x^k+x^{w-k})+x^w$, 
$\Iealog(F)=\Ieslog(F)=\langle y, x^k\rangle$, 
and thus $\taulog(C)=\teslog(C)=k$. 
\end{proposition}
\begin{proof}
Note that $C$ has two branches intersecting each other with 
multiplicity $k$. 
Since $w > 2k$, they intersect $D=\bV(y)$ 
with multiplicities $k$ and $w-k>k$. 
We may assume that those branches are 
$C_1=\bV(y+x^k)$ and $C_2=\bV(y+x^{w-k}U(x))$, 
$U(x)\in\bC[[x]]^\times$. 
By Lemma \ref{lem_scaling} we find $v(x)\in\bC[[x]]^\times$ 
such that $v(x)^{w-2k}U(xv(x))=1$. 
Substituting $xv(x)$ and $yv(x)^k$ into $x$ and $y$ 
we see that $C$ is equivalent to $\bV((y+x^k)(y+x^{w-k}))$. 

Let $f_1=\partial_x F=(w-k)x^{w-k-1}y+kx^{k-1}y+wx^{w-1}$, 
$f_2=\partial_y F=2y+x^{w-k}+x^k$ 
and $f_3=wf-x\partial_x f=wy+kx^{w-k}+(w-k)x^k$ 
where $f=y+x^{w-k}+x^k$, so that 
$\Iealog(F)=\langle f_1, f_2, f_3 \rangle$. 
Since $w\not=2k$, 
$x^{w-k}-x^k = (w-2k)^{-1}(wf_2-2f_3)$ 
is contained in $\Iealog(F)$. 
Since $-1+x^{w-2k}$ is a unit, $x^k$ is also contained. 
Then we see easily that $\Iealog(F)=\langle y, x^k\rangle$. 
Since $\dim \bC[[x, y]/\Icd(F)=k$ and 
$\Iealog(F)\subseteq \Ieslog(F)\subseteq \Icd(F)$, 
we also have $\Ieslog(F)=\langle y, x^k\rangle$. 
\end{proof}

We start with the case $\tec(C)=1$, i.e., when $C$ is nodal. 
Let $w=C.D$. 

\begin{proposition}\label{prop_node}
If $C$ is of type $A_1$ (ordinary node), 
then $C$ is equivalent to $\bV(F)$ with $F=yx+x^w$. 
In this case $\Iealog(F)=\Ieslog(F)=\langle x, y\rangle$ 
and $\teslog(C)=1$. 
\end{proposition}

\begin{proof}
If $C$ is nodal, 
then $C=C_1\cup C_2$ where $C_1.D=1$ and 
$C_2.D=w-1$. 
Let $x$ be a defining equation for $C_1$. 
Then $x, y$ form a regular parameter system on $S$, 
and $C_2$ is defined by $y=x^{w-1}u(x)$ where $u(x)\in\bC[[x, y]]^\times$. 
For $y'=-yu(x)^{-1}$, we have $C=\bV(xy'+x^w)$. 

Since $\Iea(F)=\Iec(F)=\langle x, y\rangle$ 
by Lemma \ref{lem_tec_atmost_three}, 
we also have $\Iealog(F)=\Ieslog(F)=\langle x, y\rangle$ 
by Theorem \ref{thm_inclusion}. 
Note that Proposition \ref{prop_two_branches} is also applicable if $w\geq 3$. 
\end{proof}

\begin{remark}
By Corollary \ref{cor_versal_log}, 
the semiuniversal log deformations of $\bV(yx+x^w)$ 
is given by $y(x+s)+x^w$. 
\end{remark}

\subsection{Case $\tec(C)=2$}

If $\tec(C)=2$, 
then $C$ is of type $A_2$ or $A_3$.

\begin{proposition}\label{prop_cusp}
If $C$ is of type $A_2$ (ordinary cusp), 
then $w=2$ or $w=3$. 
\begin{enumerate}[(1)]
\item
If $w=2$, $C$ is equivalent to $\bV(F)$ with $F=y^3+x^2$, 
and $\Iealog(F)=\Ieslog(F)=\langle x, y^2\rangle$
and $\teslog(C)=2$. 
\item
If $w=3$, $C$ is equivalent to $\bV(F)$ with $F=y^2+x^3$, 
and $\Iealog(F)=\Ieslog(F)=\langle x^2, y\rangle$
and $\teslog(C)=2$. 
\end{enumerate}
\end{proposition}

\begin{proof}
By considering the blow up of $S$ at $P$, 
the intersection number $w=C.D$ is seen to be $3$ or $2$ 
according as whether $D$ is tangent to $C$ or not. 
Let $E$ be a defining equation of $C$. 

\smallskip
(1)
If $w=2$, we may assume $E=x^2+xyf(y)+yg(y)$ 
by Weierstrass preparation theorem. 
By a coordinate change from $x$ to $x+(1/2)yf(y)$, 
we may assume that $E=x^2+yg(y)$
and $\mathrm{ord}_y g(y)=2$. 
By taking $y'$ to be a third root of $yg(y)$, 
we have $E=(y')^3 + x^2=F(x, y')$. 

We have 
$\Iealog(F)=\langle x, y^2\rangle$ 
by Lemma \ref{lem_binom} 
and $\Iec(F)=\langle x, y^2\rangle$ 
by Lemma \ref{lem_tec_atmost_three}. 
By Theorem \ref{thm_inclusion}, 
$\Ieslog(F)=\langle x, y^2\rangle$ 
and $\teslog(C)=2$ hold. 

\smallskip
(2)
If $w=3$, using Weierstrass preparation theorem and 
choosing $x$ appropriately, 
we may assume $E=x^3+xyf(y)+yg(y)$. 
Since $C$ is cuspidal, $E$ has no linear term and 
its quadratic term has to be the square of a nonzero linear form. 
Thus we have $\mathrm{ord}_y g(y)=1$ and $\mathrm{ord}_y f(y)\geq 1$. 
By replacing $y$ by a square root of 
$yg(y)+xyf(y)$, we may assume $E=x^3+y^2=F$. 

By Lemma \ref{lem_binom}, 
Lemma \ref{lem_tec_atmost_three} 
and Theorem \ref{thm_inclusion}, 
we have 
$\Iealog(F)=\Ieslog(F)=\Iec(F)=\langle x^2, y\rangle$
and $\teslog(C)=2$. 
\end{proof}

\begin{remark}
\begin{enumerate}[(1)]
\item
By Corollary \ref{cor_versal_log}, 
the semiuniversal log deformations of $\bV(y^3+x^2)$ and $\bV(y^2+x^3)$ 
are given by 
$y(y^2+sy+t)+x^2$ and $y(y+sx+t)+x^3$. 
\item
We can regard these families as families of plane curves over an algebraic scheme 
(cf. \cite[Corollary 3.20, Lemma 3.21, Theorem 3.27]{DiazHarris1988}). 
Then the equisingular loci are both $s=t=0$: 
the singular locus is $t=0$, and for $t=0$ and $s\not=0$, 
the curve has $2$ tangent directions at the origin. 
The statements on $\Ieslog(F)$ also follow from this. 
\end{enumerate}
\end{remark}

\begin{proposition}\label{prop_a3}
Let $C$ be of type $A_3$ (a tacnode). 
Then $w=2$ or $w\geq 4$. 
\begin{enumerate}[(1)]
\item
If $w=2$, then $C$ is equivalent to $\bV(F)$ with $F=y^4+x^2$, 
$\Iealog(F)=\Ieslog(F)=\langle x, y^3\rangle$
and $\teslog(C)=3$. 
\item
If $w=4$, then $C$ is equivalent to $\bV(F_a)$ with $F_a=y^2+ax^2y+x^4$ ($a\not=\pm 2$), 
$\Iealog(F_a)=\langle 2y+ax^2, x^3\rangle$, 
$\Ieslog(F_a)=\langle y, x^2\rangle$ 
and $\teslog(C)=2$. 
In the semiuniversal log deformation $y(y+(a+s)x^2+tx+u)+x^4$, 
the equisingular locus is $t=u=0$. 
\item
If $w\geq 5$, 
then $C$ is equivalent to $\bV(F)$ with $F=(y+x^2)(y+x^{w-2})=y(y+x^2+x^{w-2})+x^w$, 
$\Iealog(F)=\Ieslog(F)=\langle y, x^2\rangle$, 
and $\teslog(C)=2$. 
\end{enumerate}
\end{proposition}
\begin{proof}
Since the two branches $C_1, C_2$ of $C$ are tangent to each other, 
the intersection number $w=C.D$ must be $2$ or equal to or greater than $4$. 

\smallskip
(1)
If $w=2$, then the two branches are transversal to $D$. 
Let $x$ be a defining equation for $C_1$. 
Then $x, y$ form a regular parameter system 
and $C_2$ can be written as $x=y^2u(y)$ with $u(y)\in\bC[[y]]^\times$. 
If $y'=y\sqrt{u(y)}$ and $x'=\sqrt{-1}(-2x+y^2u(y))$, 
then $(y')^4+(x')^2=-4x(x-y^2u(y))$ and this defines $C$. 

By Lemma \ref{lem_binom}, 
$\Iea(F)=\Iealog(F)=\langle x, y^3\rangle$. 
The multiplicities of essential infinitely near points of $\bV(yF)$ are 
$3$ (free point) and $2$ (free point), so 
$\tes(yF)=6$ by Theorem \ref{thm_tes} 
and $\teslog(F) = 3$ by Corollary \ref{cor_tes_diff}. 
Thus $\Ieslog(F)=\Iealog(F)$. 

\smallskip
(2)
If $w=4$, 
then we take coordinates $x, y$ such that 
$C_1=\bV(y+x^2)$ and $C_2=\bV(y+x^2u(x))$, $u(0)\not=0, 1$. 
Let $v(x)$ be the square root of $u(x)/u(0)$ with $v(0)=1$ and 
\[
x'=x+\frac{v(x)-1}{x(1-u(x))} (y+x^2). 
\]
Then $x', y$ form a regular parameter system. 
Since $x'\equiv x \mod y+x^2$, we have $C_1=\bV(y+(x')^2)$. 
On the other hand, on $C_2$ we have $y=-x^2u(x)$ and 
\[
x'=x+\frac{v(x)-1}{x(1-u(x))} (-x^2u(x)+x^2)=xv(x), 
\]
hence $y+u(0)(x')^2=y+x^2u(x)=0$. 
Thus $C$ is defined by $(y+(x')^2)(y+u(0)(x')^2)$, $u(0)\not=0, 1$. 
With $y=\sqrt{u(0)}y'$ we have the equation 
$F_a(x', y')$ with $a=\sqrt{u(0)}+\sqrt{u(0)}^{-1}$. 
Note that Lemma \ref{lem_trinom} does not apply here. 

It is straightforward to see that 
$\Iea(F_a)=\Iealog(F_a)=\langle axy + 2x^3, 2y+ax^2\rangle$. 
This ideal contains $2(axy+2x^3)-ax(2y+ax^2) = (4-a^2)x^3$, 
and hence $x^3$. 
It is easy to see that $\{2y+ax^2, x^3\}$ is a Gr\"obner basis for $>_{y>x}$.  

The essential infinitely near points of $\bV(yF_a)$ are 
two free points of multiplicity $3$, 
so $\tes(yF_a)=9$ by Theorem \ref{thm_tes}, 
and $\teslog(F_a)=2$ by Corollary \ref{cor_tes_diff}. 
By Lemma \ref{lem_tec_atmost_three} 
and Theorem \ref{thm_inclusion}, 
we have 
$\Ieslog(F_a)=\Iec(F_a)=\langle y, x^2\rangle$. 

\smallskip
(3)
This is Proposition \ref{prop_two_branches} with $k=2$. 
\end{proof}

\subsection{Case $\tec(C)=3$}

If $\tec(C)=3$, 
then $C$ is of type $A_4$, $A_5$ or $D_4$. 

\begin{proposition}\label{prop_a4}
Let $C$ be of type $A_4$. 
Then $w$ is $2, 4$ or $5$. 
\begin{enumerate}[(1)]
\item
If $w=2$, then $\teslog(C)=4$. 
\item
If $w=4$, then $C$ is equivalent to $\bV(F)$ with $F=y^2+(2x^2+x^3)y+x^4$, 
$\Iealog(F)=\Ieslog(F)=\langle y+x^2, x^3\rangle$, 
and $\teslog(C)=3$. 
\item
If $w=5$, 
then $C$ is equivalent to $\bV(F_a)$ with $F_a=y^2+ax^3y+x^5$.
 
For $a=0$, $\Iealog(F_0)=\langle y, x^4\rangle$ and 
$\Ieslog(F_0)=\langle y, x^3\rangle$. 
If $a\not=0$, $\bV(F_a)$ is equivalent to $\bV(F_1)$, 
and $\Iealog(F_a)=\Ieslog(F_a)=\langle y, x^3\rangle$. 
In both cases  $\teslog(C)=3$. 
\end{enumerate}
\end{proposition}
\begin{proof}
(1)
An  embedded resolution of $C\cup D$ is given by $4$ blowups, 
and 
the multiplicities of proper transforms of $C\cup D$ 
at the essential points are 
$3$ (free point), $2$ (free point), $1$ (free point) and $1$. 
Thus $\tes(C\cup D)=7$ by Theorem \ref{thm_tes} 
and $\teslog(C)=4$ by Corollary \ref{cor_tes_diff}. 

\smallskip
(2)
If $w=4$, 
then the unique tangent direction of $C$ is equal to that of $D$, 
and we may assume that $C$ is defined by $E=y(y+f(x, y))+x^4$ 
with $f(x, y)\in\langle x, y\rangle^2$. 
Blowing up with $y=xy_1$, 
the proper transform of $C$ is given by $y_1^2 +xy_1 f_1(x, y_1)+x^2=0$ 
with $f(x, xy_1)=x^2f_1(x, y_1)$. 
This should again have a unique tangent direction, 
so $f_1(0, 0)=\pm2$. 
By replacing $y$ by $-y$ if necessary, 
we have $f=2x^2+x^3g(x)+yh(x, y)$ with $h(x, y)\in\langle x, y\rangle$ 
and $E=y(y+yh(x, y)+2x^2+x^3g(x))+x^4$. 

If we replace $y$ by $y\sqrt{1+h(x, y)}$, 
the equation is transformed to $y(y+2x^2+x^2A(x, y))+x^4$ 
with $A(x, y)\in\langle x, y\rangle$. 
We inductively show that it is equivalent to 
$y(y+2x^2+x^3A_i(x)+x^{2i}yB_i(x, y))+x^4$ 
for $i\geq 1$: 
the case $i=1$ follows from $A(x, y)\in\langle x, y\rangle$. 
Assuming the case $i\geq 1$, 
for $y'=y\sqrt{1+x^{2i}B_i(x, y)}$, 
we have 
\[
y(y+2x^2+x^3A_i(x)+x^{2i}yB_i(x, y))+x^4
= y' (y' + (1+x^{2i}B_i(x, y))^{-\frac{1}{2}}(2x^2+x^3A_i(x)))+x^4. 
\]
By writing $B_i(x, y)$ as $C_i(x)+y'D_i(x, y')$, 
we see that this can be rewritten as 
$y'(y'+2x^2+x^3A_{i+1}(x)+x^{2(i+1)}y' B_{i+1}(x, y'))+x^4$. 
Taking the limit, we may assume that $C$ is defined by $y(y+2x^2+x^3A(x))+x^4$. 
This is equal to $(y+x^2+(1/2)x^3A(x))^2-x^5A(x)-(1/4)x^6A(x)^2$, 
so $A(0)\not=0$. 
By Lemma \ref{lem_scaling}, 
we can find $v(x)\in\bC[[x]]^\times$ such that $v(x)A(xv(x))=1$. 
By substituting $xv(x)$ and $yv(x)^2$ into $x$ and $y$, 
we have the equation $y(y+2x^2+x^3)+x^4$. 

Let $f_1=\partial_x F=4x^3+3x^2y+4xy$, 
$f_2=\partial_y F=x^3+2x^2+2y$ 
and $f_3=4f-x\partial_x f=x^3+4x^2+4y$ 
where $f=y+2x^2+x^3$, so that 
$\Iealog(F)=\langle f_1, f_2, f_3\rangle$. 
From $f_2$ and $f_3$, 
we see that $y+x^2$ and $x^3$ are contained. 
Now it is easy to see that $\Iealog(F)=\langle y+x^2, x^3\rangle$. 

From $\taulog(F)=\tec(F)=3$, 
we have 
$\Iealog(F)=\Ieslog(F)=\Iec(F)$ 
and $\teslog(F)=3$. 

\smallskip
(3)
As in (2), we may assume $E=y(y+f(x, y))+x^5$ 
with $f(x, y)\in\langle x, y\rangle^2$. 
Blowing up with $y=xy_1$, 
we have $y_1^2 +xy_1 f_1(x, y_1)+x^3=0$ 
with $f(x, xy_1)=x^2f_1(x, y_1)$. 
This has to be a cusp tangent to $y_1=0$, 
so $f_1(0, 0)=0$, and 
$E=y(y+x^kg(x)+yh(x, y))+x^5$ with $h(x, y)\in\langle x, y\rangle$, 
$g(0)\not=0$ and $k\geq 3$ or $g(x)=0$. 
If $k\geq 4$, by substituting $x-(1/5)x^{k-4}yg(x)$ to $x$, 
we may assume $g(x)=0$. 
Then this is equivalent to $y^2+x^5$. 

If not, then $k=3$ and 
$E=y(y+x^3g(x)+yh(x, y))+x^5$ with $g(0)\not=0$, $h(x, y)\in\langle x, y\rangle$. 
As in (2), 
we can eliminate the $y$-part of $x^3g(x)+yh(x, y)$ 
to obtain an equation $y(y+x^3U(x))+x^5$, $U(0)\not=0$. 
By Lemma \ref{lem_scaling} we can find a unit $v(x)$ such that $v(x)U(xv(x)^2)=1$, 
and then substituting $xv(x)^2$ and $yv(x)^5$ into $x$ and $y$, 
we have $y(y+x^3)+x^5$. 

In both cases, the essential points for $\bV(yF_a)$ are of multiplicities 
$3$ (free), $3$ (free), $2$ (free) and $1$ (not free), 
so $\tes(yF_a)=12$, 
and $\teslog(F_a)=3$ by Corollary \ref{cor_tes_diff}. 

We have $\Iealog(F_0)=\langle y, x^4\rangle$ 
by Lemma \ref{lem_binom}, 
and there is a $1$-dimensional equisingular family. 
The family $y(y+sx^3)+x^5$ is equisingular, 
so $\Ieslog(F_0)=\langle y, x^3\rangle$. 
On the other hand, $\Iealog(F_1)=\langle 5x^4+3x^2y, x^3, y\rangle=\langle y, x^3\rangle$ 
by Lemma \ref{lem_trinom}. 
From $\teslog(F_a)=3$, we have $\Ieslog(F_a)=\Iealog(F_a)$ for $a\not=0$. 
\end{proof}

\begin{proposition}\label{prop_a5}
Let $C$ be of type $A_5$. 
Then $w=2, 4$ or $w\geq 6$. 
\begin{enumerate}[(1)]
\item
If $w=2$, then $\teslog(C)=5$. 
\item
If $w=4$, then $\teslog(C)=4$. 
\item
If $w=6$, 
then $C$ is equivalent to $\bV(F_a)$ or $\bV(G_a)$, $a\not=0, 1$, with 
 $F_a=(y+x^3)(y+ax^3)$ and 
$G_a=(y+x^3)(y+ax^3+x^4)$. 
In both cases  $\teslog(C)=3$. 

We have 
$\Iealog(F_a)=\langle 2y+(a+1)x^3, x^5\rangle$, 
$\Ieslog(F_a)=\langle y, x^3\rangle$, 
$\Iealog(G_a)=\langle 2y+(a+1)x^3, x^4\rangle$,
$\Ieslog(G_a)=\langle y, x^3\rangle$. 
\item
If $w\geq 7$, 
$C$ is equivalent to $\bV(F)$ with $F=(y+x^3)(y+x^{w-3})=y(y+x^3+x^{w-3})+x^w$, 
$\Iealog(F)=\Ieslog(F)=\langle y, x^3\rangle$, 
and $\teslog(C)=3$. 
\end{enumerate}
\end{proposition}
\begin{proof}
(1)
The multiplicities of proper transforms of $C\cup D$ 
at the essential points are 
$3$ (free), $2$ (free) and $2$ (free), 
so $\tes(C\cup D)=8$ by Theorem \ref{thm_tes} 
and $\teslog(C)=5$ by Corollary \ref{cor_tes_diff}. 

\smallskip
(2)
The multiplicities of proper transforms of $C\cup D$ 
are 
$3$ (free), $3$ (free) and $2$ (free), 
so $\tes(C\cup D)=11$ by Theorem \ref{thm_tes} 
and $\teslog(C)=4$ by Corollary \ref{cor_tes_diff}. 

\smallskip
(3)
It is easy to see that 
$C$ is equivalent to $\bV((y+x^3)(y+ax^3+x^4h(x)))$, $a\not=0, 1$. 

If $h(0)\not=0$, 
we consider substitutions $xv(x)$ and $yv(x)^3$ 
into $x$ and $y$ for a unit $v(x)$, 
and the new equation is 
$(y+x^3)(y+ax^3+x^4v(x)h(xv(x)))$. 
By Lemma \ref{lem_scaling}, 
we may choose $v(x)$ so that 
the equation becomes $(y+x^3)(y+ax^3+x^4)$. 

If $h(0)=0$, 
we proceed as in the proof of Proposition \ref{prop_a3} (2). 
Let $v(x)$ be the cube root of $1+xh(x)/a$ with $v(0)=1$. 
Then $v(x)\equiv 1\mod x^2$. Let 
\[
x'=x+\frac{v(x)-1}{x^2((1-a)-xh(x))} (y+x^3). 
\]
From $x'\equiv x \mod (y+x^3)$, we have $\bV(y+x^3)=\bV(y+(x')^3)$. 
Modulo $y+ax^3+x^4h(x)$, 
\[
x'\equiv x+\frac{v(x)-1}{x^2((1-a)-xh(x))}(-ax^3-x^4h(x)+x^3)=xv(x), 
\]
and $y+a(x')^3\equiv y+ax^3(1+xh(x)/a)\equiv 0$. 
Thus $C=\bV((y+(x')^3)(y+a(x')^3))$, $a\not=0, 1$. 

It is easy to see that 
$\Iealog(F_a)=\langle 2y+(a+1)x^3, x^5\rangle$. 
On the other hand, 
$\Iealog(G_a)$ is generated by $f_1:=2y+(a+1)x^3+x^4$, 
$f_2:=3(a+1)x^2y+4x^3y+6ax^5+7x^6$ 
and 
\begin{eqnarray*}
f_3 & := & 
(a+x)^2\left(6\cdot
\frac{y+(a+1)x^3+x^4}{a+x}-x\partial_x \frac{y+(a+1)x^3+x^4}{a+x}\right) \\
& = & 
(6a+7x)y+3a(a+1)x^3+(6a+4)x^4+3x^5. 
\end{eqnarray*}
Then $\Iealog(G_a)$ contains 
$(6a+7x)f_1-2f_3=(a-1)x^4+x^5$, 
and hence $x^4$. 
Now it is easy to see that
$\Iealog(G_a)=\langle 2y+(a+1)x^3, x^4\rangle$. 

From Lemma \ref{lem_tec_atmost_three} 
we deduce that 
$\Iec(F_a)=\Iec(G_a)=\langle y, x^3\rangle$ 
using coordinate changes that map $\bV(y^2+x^6)$ to $\bV(F_a)$ and $\bV(G_a)$, 
or by repeating the arguments of the proof of Lemma \ref{lem_tec_atmost_three}. 

In both cases, 
the multiplicities of proper transforms of $C\cup D$ 
are 
$3$ (free), $3$ (free) and $3$ (free), 
so $\tes(C\cup D)=14$ by Theorem \ref{thm_tes} 
and $\teslog(C)=3$ by Corollary \ref{cor_tes_diff}. 
Thus $\Ieslog(F_a)=\Iec(F_a)$ 
and $\Ieslog(G_a)=\Iec(G_a)$ hold.

\smallskip
(4)
This is Proposition \ref{prop_two_branches} with $k=3$. 
\end{proof}

\begin{proposition}\label{prop_d4}
If $C$ is of type $D_4$, then $w\geq 3$. 
\begin{enumerate}[(1)]
\item
If $w=3$, then 
$C$ is equivalent to $\bV(F_a)$ with $F_a=x(y+x)(y+ax)$, $a\not=0, 1$, 
and 
$\Iealog(F_a)=
\langle y^2-(a^2-a+1)x^2, 2xy+(a+1)x^2, x^3\rangle$, 
$\Ieslog(F_a)=\langle x^2, xy, y^2\rangle$, 
and $\teslog(C)=3$. 
\item
If $w\geq 4$, then 
$C$ is equivalent to $\bV(F)$ with $F=x(y+x)(y+x^{w-2})=y(yx+x^2+x^{w-1})+x^w$, 
$\Iealog(F)=\Ieslog(F)=\langle x^2, xy, y^2\rangle$, 
and $\teslog(C)=3$. 
\end{enumerate}
\end{proposition}
\begin{proof}
Since the multiplicity of $C$ is $3$, 
we have $w\geq 3$. 

\smallskip
(1)
The union $C\cup D$ is an ordinary quadruple point, 
and hence is isomorphic to the completion of 
the union of $4$ distinct lines by \cite[p.\ 260, {\bf 14}]{AGV1988}. 
Thus $C$ is equivalent to $\bV(x(y+x)(y+ax))$, $a\not=0, 1$. 

It is easy to see that 
$\Iealog(F_a)$ is generated by 
$f_1:=y^2+2(a+1)xy+3ax^2$ and $f_2:=2xy+(a+1)x^2$. 
Then 
$\Iealog(F_a)$ contains 
$(2y+3(a+1)x)f_2-4xf_1=3(a-1)^2x^3$, 
and hence $f_3:=x^3$. 
For the order $>_{y>x}$, 
$\{f_1, f_2, f_3\}$ is a Gr\"obner basis. 
We can further replace $f_1$ by $f_1-(a+1)f_2=y^2-(a^2-a+1)x^2$. 

We obtain an embedded resolution by 
blowing up a (free) point of multiplicity $4$, 
so $\tes(C\cup D)=8$ by Theorem \ref{thm_tes} 
and $\teslog(C)=3$ by Corollary \ref{cor_tes_diff}. 
Since $\tec(C)=3$, 
$\Ieslog(F_a)=\Iec(F_a)=\langle x^2, xy, y^2\rangle$ holds. 

\smallskip
(2)
We may assume that the $3$ branches of $C$ are 
$\bV(x)$, $\bV(y+x)$ and $\bV(y+x^{w-2}U(x))$ 
with $U(x)\in\bC[[x]]^\times$. 
Taking $v(x)\in\bC[[x]]^\times$ 
such that $v(x)^{w-3}U(xv(x))=1$ by Lemma \ref{lem_scaling} 
and substituting $xv(x)$ and $yv(x)$ into $x$ and $y$, 
we have $C=\bV(x(y+x)(y+x^{w-2}))$. 

We can write $F=yf+x^w$ with $f=x^{w-1}+x^2+xy$, 
and $\Iealog(F)$ is generated by 
$f_1:=\partial_x F = wx^{w-1}+(w-1)x^{w-2}y+2xy+y^2$, 
$f_2:=\partial_y F = x^{w-1}+x^2+2xy$ 
and 
$f_3:= wf-x\partial_x f= x^{w-1}+(w-2)x^2+(w-1)xy$. 
From $2f_3-(w-1)f_2=(w-3)(x^2-x^{w-1})$, 
we see that $\Iealog(F)$ contains $x^2-x^{w-1}$, and hence $x^2$. 
Reducing $f_2$ modulo $x^2$ we have $xy\in \Iealog(F)$, 
and reducing $f_1$ modulo $\langle x^2, xy\rangle$ we have $y^2\in \Iealog(F)$. 
It is easy to see that 
$\Iealog(F)=\langle x^2, xy, y^2\rangle$. 

Since $\tec(C)=3$, 
$\Ieslog(F)=\langle x^2, xy, y^2\rangle$ 
and $\teslog(C)=3$ also hold. 
\end{proof}

From these results, 
singularities $C$ with $\teslogd{D}(C)\leq 3$ 
are classified as in Table \ref{table_classification}. 
In the $A_4$ case with $w=5$ and 
$A_5$ case with $w=6$, 
the ones in the upper rows 
appear as equisingular limits of the ones in the lower row. 

\begin{table}[htbp]
\centering
\caption{Singularities $C$ with $\teslogd{D}(C)\leq 3$}
\label{table_classification}
\begin{tabular}{|c|c|c|c|c|l|l|} \hline
   $\teslogd{D}$ & $\tec$ & $\delta$ & Type & $w$ & Equation & Details \\ 
\hline \hline
   $1$ & $1$ & $1$ & $A_1$ & $\geq 2$ & $yx+x^w$ & \ref{prop_node} \\ 
\hline
   $2$ & $2$ & $1$ & $A_2$ & $2$ & $y^3+x^2$ & \ref{prop_cusp} \\ 
\cline{5-6}
         &         &        &            & $3$ & $y^2+x^3$ &  \\ 
\cline{3-7}
         &         & $2$ &  $A_3$ & $4$ & $y^2+ax^2y+x^4$ ($a\not=\pm 2$) 
                   & \ref{prop_a3} (2) \\ 
\cline{5-7}
         &        &         &             & $\geq 5$ & $(y+x^2)(y+x^{w-2})$ & \ref{prop_a3} (3) \\ 
\hline
   $3$ & $2$ & $2$ & $A_3$ & $2$ & $y^4+x^2$ & \ref{prop_a3} (1) \\ 
\cline{2-7}
         & $3$ & $2$ & $A_4$ & $4$ & $y^2+(2x^2+x^3)y + x^4$ & \ref{prop_a4} (2) \\ 
\cline{5-7}
         &        &        &             & $5$ & $y^2+x^5$ & \ref{prop_a4} (3) \\ 
\cline{6-6}
         &        &        &             &        & $y^2+x^3y+x^5$ & \\ 
\cline{3-7}
         &        & $3$ & $A_5$ & $6$ & $(y+x^3)(y+ax^3)$  ($a\not=0, 1$) & \ref{prop_a5} (3) \\ 
\cline{6-6}
         &        &         &            &       & $(y+x^3)(y+ax^3+x^4)$  ($a\not=0, 1$) & \\ 
\cline{5-7}
         &        &        &            & $\geq 7$ & $(y+x^3)(y+x^{w-3})$ & \ref{prop_a5} (4) \\ 
\cline{4-7}
         &        &         & $D_4$ & $3$ & $x(y+x)(y+ax)$ ($a\not=0, 1$) & \ref{prop_d4} (1) \\ 
\cline{5-7}
         &        &         &             & $\geq 4$ & $x(y+x)(y+x^{w-2})$ & \ref{prop_d4} (2) \\ 
\hline
\end{tabular}
\end{table}

\bigskip
To conclude, let us look at a global example. 
\begin{example}
Let $B\subset \bP^2$ be a general cubic, 
$P\in B$ a point of order $4$ or $12$ for the addition on $B$ 
with an inflexion point as the unit element. 
Let $\Lambda$ be the $3$-dimensional parameter space of quartics $C$ 
such that $C|_B=12P$. 

From the table, it is expected that 
$\Lambda$ contains finitely many curves with an $A_5$- or $D_4$-singularity 
at $P$. 
Actually, according to \cite[Corollary 4.3]{Takahashi1996}, 
there is exactly $1$ member of $\Lambda$ of each type. 
\end{example}

\end{document}